\catcode`\@=11

\font\tenmsa=msam10
\font\sevenmsa=msam7
\font\fivemsa=msam5
\font\tenmsb=msbm10
\font\sevenmsb=msbm7
\font\fivemsb=msbm5
\newfam\msafam
\newfam\msbfam
\textfont\msafam=\tenmsa  \scriptfont\msafam=\sevenmsa
  \scriptscriptfont\msafam=\fivemsa
\textfont\msbfam=\tenmsb  \scriptfont\msbfam=\sevenmsb
  \scriptscriptfont\msbfam=\fivemsb

\def\hexnumber@#1{\ifnum#1<10 \number#1\else
 \ifnum#1=10 A\else\ifnum#1=11 B\else\ifnum#1=12 C\else
 \ifnum#1=13 D\else\ifnum#1=14 E\else\ifnum#1=15 F\fi\fi\fi\fi\fi\fi\fi}

\def\msa@{\hexnumber@\msafam}
\def\msb@{\hexnumber@\msbfam}
\mathchardef\boxdot="2\msa@00
\mathchardef\boxplus="2\msa@01
\mathchardef\boxtimes="2\msa@02
\mathchardef\square="0\msa@03
\mathchardef\blacksquare="0\msa@04
\mathchardef\centerdot="2\msa@05
\mathchardef\lozenge="0\msa@06
\mathchardef\blacklozenge="0\msa@07
\mathchardef\circlearrowright="3\msa@08
\mathchardef\circlearrowleft="3\msa@09
\mathchardef\rightleftharpoons="3\msa@0A
\mathchardef\leftrightharpoons="3\msa@0B
\mathchardef\boxminus="2\msa@0C
\mathchardef\Vdash="3\msa@0D
\mathchardef\Vvdash="3\msa@0E
\mathchardef\vDash="3\msa@0F
\mathchardef\twoheadrightarrow="3\msa@10
\mathchardef\twoheadleftarrow="3\msa@11
\mathchardef\leftleftarrows="3\msa@12
\mathchardef\rightrightarrows="3\msa@13
\mathchardef\upuparrows="3\msa@14
\mathchardef\downdownarrows="3\msa@15
\mathchardef\upharpoonright="3\msa@16

\mathchardef\downharpoonright="3\msa@17
\mathchardef\upharpoonleft="3\msa@18
\mathchardef\downharpoonleft="3\msa@19
\mathchardef\rightarrowtail="3\msa@1A
\mathchardef\leftarrowtail="3\msa@1B
\mathchardef\leftrightarrows="3\msa@1C
\mathchardef\rightleftarrows="3\msa@1D
\mathchardef\Lsh="3\msa@1E
\mathchardef\Rsh="3\msa@1F
\mathchardef\rightsquigarrow="3\msa@20
\mathchardef\leftrightsquigarrow="3\msa@21
\mathchardef\looparrowleft="3\msa@22
\mathchardef\looparrowright="3\msa@23
\mathchardef\circeq="3\msa@24
\mathchardef\succsim="3\msa@25
\mathchardef\gtrsim="3\msa@26
\mathchardef\gtrapprox="3\msa@27
\mathchardef\multimap="3\msa@28
\mathchardef\therefore="3\msa@29
\mathchardef\because="3\msa@2A
\mathchardef\doteqdot="3\msa@2B

\mathchardef\triangleq="3\msa@2C
\mathchardef\precsim="3\msa@2D
\mathchardef\lesssim="3\msa@2E
\mathchardef\lessapprox="3\msa@2F
\mathchardef\eqslantless="3\msa@30
\mathchardef\eqslantgtr="3\msa@31
\mathchardef\curlyeqprec="3\msa@32
\mathchardef\curlyeqsucc="3\msa@33
\mathchardef\preccurlyeq="3\msa@34
\mathchardef\leqq="3\msa@35
\mathchardef\leqslant="3\msa@36
\mathchardef\lessgtr="3\msa@37
\mathchardef\backprime="0\msa@38
\mathchardef\risingdotseq="3\msa@3A
\mathchardef\fallingdotseq="3\msa@3B
\mathchardef\succcurlyeq="3\msa@3C
\mathchardef\geqq="3\msa@3D
\mathchardef\geqslant="3\msa@3E
\mathchardef\gtrless="3\msa@3F
\mathchardef\sqsubset="3\msa@40
\mathchardef\sqsupset="3\msa@41
\mathchardef\trianglerighteq="3\msa@44
\mathchardef\trianglelefteq="3\msa@45
\mathchardef\bigstar="0\msa@46
\mathchardef\between="3\msa@47
\mathchardef\blacktriangledown="0\msa@48
\mathchardef\blacktriangleright="3\msa@49
\mathchardef\blacktriangleleft="3\msa@4A
\mathchardef\blacktriangle="0\msa@4E
\mathchardef\triangledown="0\msa@4F
\mathchardef\eqcirc="3\msa@50
\mathchardef\lesseqgtr="3\msa@51
\mathchardef\gtreqless="3\msa@52
\mathchardef\lesseqqgtr="3\msa@53
\mathchardef\gtreqqless="3\msa@54
\mathchardef\Rrightarrow="3\msa@56
\mathchardef\Lleftarrow="3\msa@57
\mathchardef\veebar="2\msa@59
\mathchardef\barwedge="2\msa@5A
\mathchardef\doublebarwedge="2\msa@5B
\mathchardef\angle="0\msa@5C
\mathchardef\measuredangle="0\msa@5D
\mathchardef\sphericalangle="0\msa@5E
\mathchardef\varpropto="3\msa@5F
\mathchardef\smallsmile="3\msa@60
\mathchardef\smallfrown="3\msa@61
\mathchardef\Subset="3\msa@62
\mathchardef\Supset="3\msa@63
\mathchardef\Cup="2\msa@64

\mathchardef\Cap="2\msa@65

\mathchardef\curlywedge="2\msa@66
\mathchardef\curlyvee="2\msa@67
\mathchardef\leftthreetimes="2\msa@68
\mathchardef\rightthreetimes="2\msa@69
\mathchardef\subseteqq="3\msa@6A
\mathchardef\supseteqq="3\msa@6B
\mathchardef\bumpeq="3\msa@6C
\mathchardef\Bumpeq="3\msa@6D
\mathchardef\lll="3\msa@6E

\mathchardef\ggg="3\msa@6F

\mathchardef\circledS="0\msa@73
\mathchardef\pitchfork="3\msa@74
\mathchardef\dotplus="2\msa@75
\mathchardef\backsim="3\msa@76
\mathchardef\backsimeq="3\msa@77
\mathchardef\complement="0\msa@7B
\mathchardef\intercal="2\msa@7C
\mathchardef\circledcirc="2\msa@7D
\mathchardef\circledast="2\msa@7E
\mathchardef\circleddash="2\msa@7F
\def\ulcorner{\delimiter"4\msa@70\msa@70 }
\def\urcorner{\delimiter"5\msa@71\msa@71 }
\def\llcorner{\delimiter"4\msa@78\msa@78 }
\def\lrcorner{\delimiter"5\msa@79\msa@79 }
\def\yen{\mathhexbox\msa@55 }
\def\checkmark{\mathhexbox\msa@58 }
\def\circledR{\mathhexbox\msa@72 }
\def\maltese{\mathhexbox\msa@7A }
\mathchardef\lvertneqq="3\msb@00
\mathchardef\gvertneqq="3\msb@01
\mathchardef\nleq="3\msb@02
\mathchardef\ngeq="3\msb@03
\mathchardef\nless="3\msb@04
\mathchardef\ngtr="3\msb@05
\mathchardef\nprec="3\msb@06
\mathchardef\nsucc="3\msb@07
\mathchardef\lneqq="3\msb@08
\mathchardef\gneqq="3\msb@09
\mathchardef\nleqslant="3\msb@0A
\mathchardef\ngeqslant="3\msb@0B
\mathchardef\lneq="3\msb@0C
\mathchardef\gneq="3\msb@0D
\mathchardef\npreceq="3\msb@0E
\mathchardef\nsucceq="3\msb@0F
\mathchardef\precnsim="3\msb@10
\mathchardef\succnsim="3\msb@11
\mathchardef\lnsim="3\msb@12
\mathchardef\gnsim="3\msb@13
\mathchardef\nleqq="3\msb@14
\mathchardef\ngeqq="3\msb@15
\mathchardef\precneqq="3\msb@16
\mathchardef\succneqq="3\msb@17
\mathchardef\precnapprox="3\msb@18
\mathchardef\succnapprox="3\msb@19
\mathchardef\lnapprox="3\msb@1A
\mathchardef\gnapprox="3\msb@1B
\mathchardef\nsim="3\msb@1C
\mathchardef\napprox="3\msb@1D
\mathchardef\nsubseteqq="3\msb@22
\mathchardef\nsupseteqq="3\msb@23
\mathchardef\subsetneqq="3\msb@24
\mathchardef\supsetneqq="3\msb@25
\mathchardef\subsetneq="3\msb@28
\mathchardef\supsetneq="3\msb@29
\mathchardef\nsubseteq="3\msb@2A
\mathchardef\nsupseteq="3\msb@2B
\mathchardef\nparallel="3\msb@2C
\mathchardef\nmid="3\msb@2D
\mathchardef\nshortmid="3\msb@2E
\mathchardef\nshortparallel="3\msb@2F
\mathchardef\nvdash="3\msb@30
\mathchardef\nVdash="3\msb@31
\mathchardef\nvDash="3\msb@32
\mathchardef\nVDash="3\msb@33
\mathchardef\ntrianglerighteq="3\msb@34
\mathchardef\ntrianglelefteq="3\msb@35
\mathchardef\ntriangleleft="3\msb@36
\mathchardef\ntriangleright="3\msb@37
\mathchardef\nleftarrow="3\msb@38
\mathchardef\nrightarrow="3\msb@39
\mathchardef\nLeftarrow="3\msb@3A
\mathchardef\nRightarrow="3\msb@3B
\mathchardef\nLeftrightarrow="3\msb@3C
\mathchardef\nleftrightarrow="3\msb@3D
\mathchardef\divideontimes="2\msb@3E
\mathchardef\varnothing="0\msb@3F
\mathchardef\nexists="0\msb@40
\mathchardef\mho="0\msb@66
\mathchardef\thorn="0\msb@67
\mathchardef\beth="0\msb@69
\mathchardef\gimel="0\msb@6A
\mathchardef\daleth="0\msb@6B
\mathchardef\lessdot="3\msb@6C
\mathchardef\gtrdot="3\msb@6D
\mathchardef\ltimes="2\msb@6E
\mathchardef\rtimes="2\msb@6F
\mathchardef\shortmid="3\msb@70
\mathchardef\shortparallel="3\msb@71
\mathchardef\smallsetminus="2\msb@72
\mathchardef\thicksim="3\msb@73
\mathchardef\thickapprox="3\msb@74
\mathchardef\approxeq="3\msb@75
\mathchardef\succapprox="3\msb@76
\mathchardef\precapprox="3\msb@77
\mathchardef\curvearrowleft="3\msb@78
\mathchardef\curvearrowright="3\msb@79
\mathchardef\digamma="0\msb@7A
\mathchardef\varkappa="0\msb@7B
\mathchardef\hslash="0\msb@7D
\mathchardef\hbar="0\msb@7E
\mathchardef\backepsilon="3\msb@7F
\def\Bbb{\ifmmode\let\next\Bbb@\else
 \def\next{\errmessage{Use \string\Bbb\space only in math mode}}\fi\next}
\def\Bbb@#1{{\Bbb@@{#1}}}
\def\Bbb@@#1{\fam\msbfam#1}

\catcode`\@=\active

\def\inv{^{\raise.15ex\hbox{${
  \scriptscriptstyle -}$}\kern-.05em 1}}

\def\Dsl{\,\raise.15ex\hbox{$/$}\mkern-13.5mu D}
\def\dsl{\raise.15ex\hbox{$/$}\kern-.57em\hbox{$\partial$}}

\def\lspace{\ifx\answ\bigans{}\else\qquad\fi}

\def\CA{\hbox{{$\cal A$}}} 
 \def\CG{\hbox{{$\cal G$}}}

\def\lform{\hbox{$\sqcup$}\llap{\hbox{$\sqcap$}}}
\def\darr#1{\raise1.5ex\hbox{$\leftrightarrow$}
\mkern-16.5mu #1}
\def\INT{{\textstyle \int\kern-.642em\int}}

\def\Z{{\Bbb Z}}

\def\eps{{\epsilon}}

\def\lcross{{>\!\!\!\triangleleft}}

\def\rbiprod{{\cdot\kern-.33em\triangleright\!\!\!<}}
\def\lbiprod{{>\!\!\!\triangleleft\kern-.33em\cdot}}

\def\tens{\mathop{\otimes}}

\def\la{{\triangleright}}
\def\cora{{\blacktriangleleft}}

\def\isom{{\cong}}

\def\span{{\rm span}}

\def\Ad{{\rm Ad}}

\def\id{{\rm id}}

\def\<{\langle}
\def\>{\rangle}

\def\vecs{{\bf s}}

\def\<{\langle}
\def\>{\rangle}

\def\equad{\kern -1.7em}
\def\qqquad{\qquad\quad}

\def\und#1{{\underline {#1}}}

\def\text#1{\mbox{\rm #1}}
\def\note#1{}

\def\blacksquare{{\lform}}%AMS Tex Fakes
\def\frac#1#2{{{#1\over#2}}}

\def\proof{\goodbreak\noindent{\bf Proof\quad}}

\def\endproof{{\ $\lform$}\bigskip }

\def\align#1{\begin{eqnarray*}#1\end{eqnarray*}}
\documentstyle[11pt, epsf]{article}
\newtheorem{lemma}{Lemma}[section] \newtheorem{propos}[lemma]{Proposition}
\newtheorem{example}[lemma]{Example} 
 \newtheorem{corol}[lemma]{Corollary}

\begin{document}\baselineskip 22pt

{\ }\qquad\qquad \hskip 4.3in  Damtp/96-31
\vspace{.2in}

\begin{center} {\LARGE  DIAGRAMMATICS OF BRAIDED GROUP GAUGE THEORY}
\\ \baselineskip 13pt{\ }
{\ }\\
 S. Majid\footnote{Royal Society University Research Fellow and Fellow of
Pembroke College, Cambridge}\\
{\ }\\
Department of Mathematics, Harvard University\\
Science Center, Cambridge MA 02138, USA\footnote{During 1995+1996}\\
+\\
Department of Applied Mathematics \& Theoretical Physics\\
University of Cambridge, Cambridge CB3 9EW\\
\end{center}
\begin{center}
March 1996
\end{center}

\vspace{10pt}
\begin{quote}\baselineskip 13pt
\noindent{\bf Abstract}  We develop a gauge theory or theory of bundles and connections on them at the level of braids and tangles. Extending recent algebraic work, we provide now a fully diagrammatic treatment of principal bundles,  a theory of global gauge transformations, associated braided fiber bundles and covariant derivatives on them. We describe the local structure for a concrete $\Z_3$-graded or `anyonic' realization of the theory.

\bigskip


\end{quote}
\baselineskip 22.5pt

\section{Introduction} 

There has been a lot of interest in recent years in developing some form of `noncommutative algebraic
geometry'.  Some years ago we introduced a  `braided approach' in which one keeps more closely the classical form of geometrical constructions but make them within a  braided category.  The novel aspect of this `braided mathematics'  is that algebraic information `flows' along braid and tangle diagrams much as information flows along the wiring in a computer, except that the under and over crossings are non-trivial (and generally distinct) operators. Constructions which work universally indeed take place in the braided category of braid and tangle diagrams. 

We have already shown in \cite{Ma:exa}\cite{Ma:bra}\cite{Ma:bg}\cite{Ma:tra}\cite{Ma:bos}\cite{Ma:lin}\cite{Ma:lie}\cite{Ma:sol} and several other papers the existence of group-like objects or `braided groups' at this level of braids and tangles, and developed their basic theory and many applications. We refer to \cite{Ma:introp}\cite{Ma:introm} for reviews. See also the last chapter of the text \cite{Ma:book} for the concrete application of this machinery to $q$-deforming physics.

In this paper we provide a systematic treatment of `gauge theory' or the theory of bundles and connections on them in this same setting. The structure group of the bundle will be a braided group as above (a Hopf algebra in a braided category). All geometrical spaces are handled through their co-ordinate rings directly, which we allow to be arbitrary algebras in our braided category. Such a {\em braided group gauge theory} has been initiated recently in our algebraic work \cite{BrzMa:coa} with T.~Brzezinski, as an example of a more general `coalgebra gauge theory'. This coalgebra gauge theory  generalised our earlier work on quantum group gauge theory\cite{BrzMa:gau} to the level of  fibrations based essentially on algebra factorisations, with braided group gauge theory as an example. We develop now purely diagrammatic proofs, thereby lifting the construction of principal bundles to a general braided category with suitable direct sums, kernels and cokernels. We also extend the theory considerably, providing now global gauge transformations, a precise characterisation of which connections come from the base in a trivial bundle, associated braided fiber bundles and the covariant derivative on their sections. This provides a fairly complete formalism  of diagrammatic braided group gauge theory, as well as a first step to developing the same ideas for the more general coalgebra gauge theory in a continuation of \cite{BrzMa:coa}. Finally, we conclude with  the local picture for a simple truly braided example of based on  $\Z_3$-graded vector spaces.

\subsection*{Acknowledgements}

I would like to thank T. Brzezinski for our continuing discussions under EPSRC research grant GR/K02244, of which this work is an offshoot.  

\subsection*{Preliminaries} For braided categories we use the conventions and notation in \cite{Ma:book}. Informally, a braided category is a collection of objects $V,W,Z,\cdots$, with a tensor product between any two which is associative up to isomorphism and commutative up to isomorphism. We omit the former isomorphism and denote the latter by $\Psi_{V,W}:V\tens W\to W\tens V$.
There are coherence theorems which ensure that these isomorphisms are consistent in the expected way. The main difference with usual or super vector spaces is that the `generalised transposition' or {\em braiding} $\Psi$ need not obey $\Psi^2=\id$. There is also a unit object $\und 1$ for the tensor product, with associated morphisms. Braided categories have been formally introduced into category theory by Joyal and Street \cite{JoyStr:bra} and also arise in the representation theory of quantum groups due  essentially to the work of Drinfeld \cite{Dri:qua}. 

An algebra in a braided category means an object $P$ and morphisms $\eta:\und 1\to P$, $\cdot:P\tens P\to P$ obeying the obvious associativity and unity axioms. A braided group means a bialgebra or Hopf algebra in a braided category, i.e. an algebra $B$ in the category and morphisms $\eps:B\to \und 1$, $\Delta:B\to B\und\tens B$ obeying the arrow-reversed algebra axioms (a coalgebra). We require that $\eps,\Delta$ are homomorphisms of braided algebras, where $B\und \tens B$ is the braided tensor product algebra. For a full braided group we usually rerquire also an antipode morphism $S:B\to B$ defined as the convolution-inverse of the identity morphism. Not only these axioms   (they are obvious enough) but the existence and construction of examples was introduced in \cite{Ma:bra}\cite{Ma:bg}.

\begin{figure}
 \[ \epsfbox{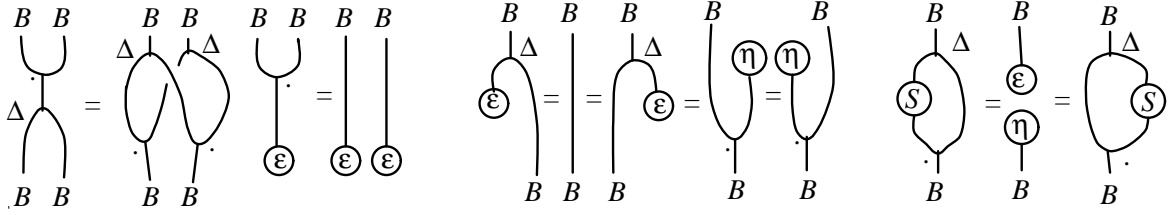}\]
 \caption{Axioms of a braided group in diagrammatic notation}
\end{figure}

The axioms of a braided group are summarized in Figure~2 in a diagrammatic notation $\Psi=\epsfbox{braid.eps}, \Psi^{-1}=\epsfbox{braidinv.eps}$ and $\cdot=\epsfbox{prodfrag.eps}$, $\Delta=\epsfbox{deltafrag.eps}$. Other morphisms are written as nodes, and the unit object $\und 1$ is denoted by omission. The functoriality of the braiding says that we can pull nodes through braid crossings as if they are beads on a string. A coherence theorem\cite{JoyStr:bra} for braided categories ensures that this notation is consistent. This technique for working with braided algebras and braided groups appeared in the 1990 work of the author \cite{Ma:bra}, and is a conjunction (for the first time) of the usual ideas of wiring diagrams in computer science (where crossings or wires have no significance) with the coherence theorem for braided categories needed for nontrivial $\Psi$. The diagrammatic theory of braided groups, actions on braided algebras, cross products by them, dual braided groups, opposite coproducts, braided-(co)commutativity, braided tensor product of representations, braided adjoint (co)action, braided Lie-algebra objects etc  have all been developed (by the author), see \cite{Ma:bra}\cite{Ma:bg}\cite{Ma:tra}\cite{Ma:bos}\cite{Ma:lie}.

In particular, we need the concept of a braided comodule algebra $P$ under a braided group $B$\cite{Ma:introp}. This means that $P$ is an algebra equipped with a morphism $\cora:P\to P\und\tens B$ which forms a comodule and which is an algebra homomorphism to the braided tensor product algebra. The diagram for the latter condition is the same as for the coproduct in Figure~1, with $\Delta$ replaced by $\cora$. 

We also assume that our category has equalisers and coequalisers compatible with the tensor product. Then associated to any comodule $P$ is a  `maximal subobject' or equaliser $P^B$ such that the `restrictions' of $\cora,\id\tens\eta:P\to P\tens B$ coincide. This means an object $P^B$ and morphism $P^B\to P$ universal with the property that the two composites $P^B\to P{ \to \atop \to} P\tens B$ coincide. Universal means that any other such  $X'\to P$ factors through $P^B$. If $P$ is a braided comodule algebra then $P^B$ is an algebra. 

If $M$ is an algebra and $P$ an $M$-bimodule we will also need  a `quotient object'   or coequaliser $P\tens_M P$ such that the two product morphisms $(\cdot_P\tens\id),(\id\tens\cdot_P):P\tens M\tens P\to P\tens P$ project to the same in $P\tens_M P$. This means an object $P\tens_M P$ and morphism $P\tens P\to P\tens_M P$ universal with the property that the two composites $P\tens M\tens P{\to\atop\to}P\tens P\to P\tens_M P$ coincide. Universal means that any other  such  $P\tens P\to X'$ factors through $P\tens_M P$. 

Most of the constructions in the paper hold at this general level by working with equalisers and coequalisers. However, for a theory of connections and differential forms, we will in fact assume that our category has suitable direct sums etc as well, i.e. an Abelian braided category. Moreover, in our diagrammatic proofs we will generally suppress the `inclusion' and `projection' morphism associated with kernels and quotients or equalisers and coequalisers, but they should be always understood where needed. Alternatively, the reader can consider that all constructions take place in a concrete braided category such as the representations of a strict quantum group.

Finally, we define differential forms on an algebra $P$ in exactly the same way as familiar in non-commutative geometry. Thus, we require an object $\Omega^1(P)$ in the braided category and morphisms $\cdot_L,\cdot_R$ by which it becomes a $P$-bimodule (suppressing throughout the implicit associativity), and a morphism $d:P\to \Omega^1(P)$ such that the Leibniz rule
\[ d\circ\cdot=\cdot_L(\id\tens d) + \cdot_R (d\tens\id).\]
holds. Higher forms are defined by the tensor product of $\Omega^1(P)$ over $P$ and the extension of $d$ by $d^2=0$.  We concentrate throughout on  $\Omega^1(P)=\Omega^1P$, the universal first order differential calculus defined as the kernel of the product morphism $\cdot_P:P\tens P\to P$ as a $P$-bimodule and $d=(\eta_P\tens\id)-(\id\tens\eta_P)$. The higher forms in this case can be identified, as in \cite{BrzMa:gau}, with the joint kernel of the $n$ morphisms $P^{\tens n+1}\to P^{\tens n}$ defined by multiplying adjacent copies of $P$.
\begin{figure}
\[ \epsfbox{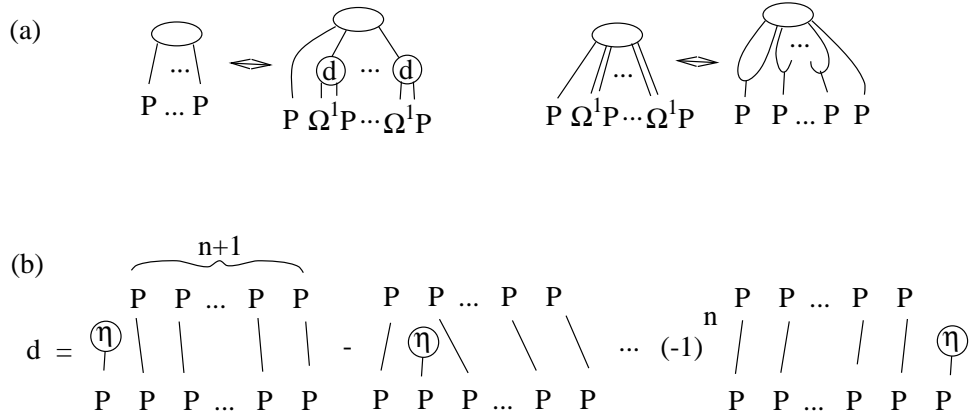}\]
\caption{(a) correspondence $\Omega^nP$ in terms of $P^{\tens n}$ and tensor product $P(\Omega^1P)^n$ (b) differential $d$ in terms of $P^{\tens n}$}
\end{figure}
We work with $\Omega^nP$ in terms of $P^{\tens n+1}$. The correspondence and the resulting differential $d$ in these terms is shown in Figure~2. In one direction we apply $d:P\to \Omega^1P$ and in the other direction we  realise each $\Omega^1P$ in $P^{\tens 2}$ and multiply up. The wedge product $\Omega^nP\tens \Omega^mP\to \Omega^{n+m}P$ is given in terms of $P^{\tens n+1}\tens P^{\tens m+1}\to P^{\tens n+m+1}$ by multiplication of the rightmost  factor in $P^{\tens n+1}$ with the leftmost factor in $P^{\tens m+1}$.

\section{General Braided Principal Bundles and Connections}

Let $B$ be a braided group in a braided category as explained above. Following \cite{BrzMa:coa}, we define a braided principal bundle with structure group $P$ to be: 

(P1) An algebra $P$ in the braided category which is a braided right $B$-comodule algebra under a right coaction $\cora$ of $B$. We denote by $M$ its associated `fixed point subalgebra' $P^B$ and suppose that $P$ is flat as
an $M$-module. 

(P2) An inverse to the morphism $\chi:P\tens_M P\to P\tens B$  inherited from $\tilde\chi=(\cdot_P\tens \id)\circ(\id\tens \circ\cora):P\tens P\to P\tens B$. The braided-comodule algebra property ensures here that $\tilde\chi$ descends to $P\tens_MP$.

Next we consider connection forms on braided principal bundles with the universal differential calculus $\Omega^1P$. Following\cite{BrzMa:coa}, this is a morphism $\omega: B\to \Omega^1P$ such that 

(C1) $\tilde\chi\circ\omega=\eta_P\tens(\id-\eta\circ\eps)$ and $\omega\circ\eta=0$

(C2) $\cora\circ\omega=(\omega\tens\id)\circ\Ad$,  where $\cora$ is the braided tensor product coaction\cite{Ma:bg} on $P\tens P$ restricted to $\Omega^1P$, and $\Ad$ is the braided adjoint coaction from \cite{Ma:lin}.

This is equivalent to an abstract definition as an equivariant splitting of $\Omega^1P$ into `horizontal forms' $P(\Omega^1M)P$ and their complement:

\begin{propos} Connections $\omega$ on a braided principal bundle $P$ are in 1-1 correspondence with morphisms $\Pi:\Omega^1P\to \Omega^1P$ which are idempotent, zero on $P(\Omega^1M)P$, left $P$-module morphisms, intertwiners for the coaction of $B$ on $\Omega^1P$ and obey $\tilde\chi\circ\Pi=\tilde\chi$, via
\[ \epsfbox{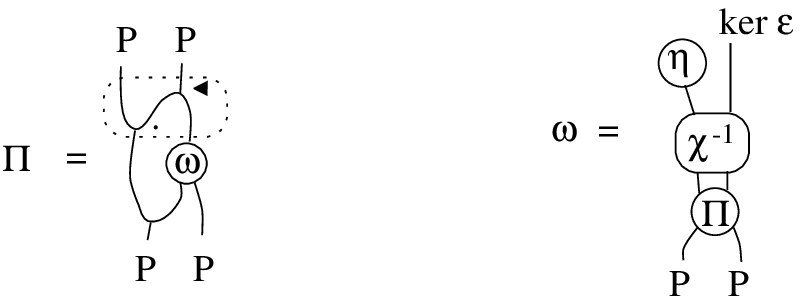}\]
and $\omega\circ\eta=0$. The box on the left is the morphism $\tilde\chi$. 
\end{propos}
\begin{figure}
 \[  \epsfbox{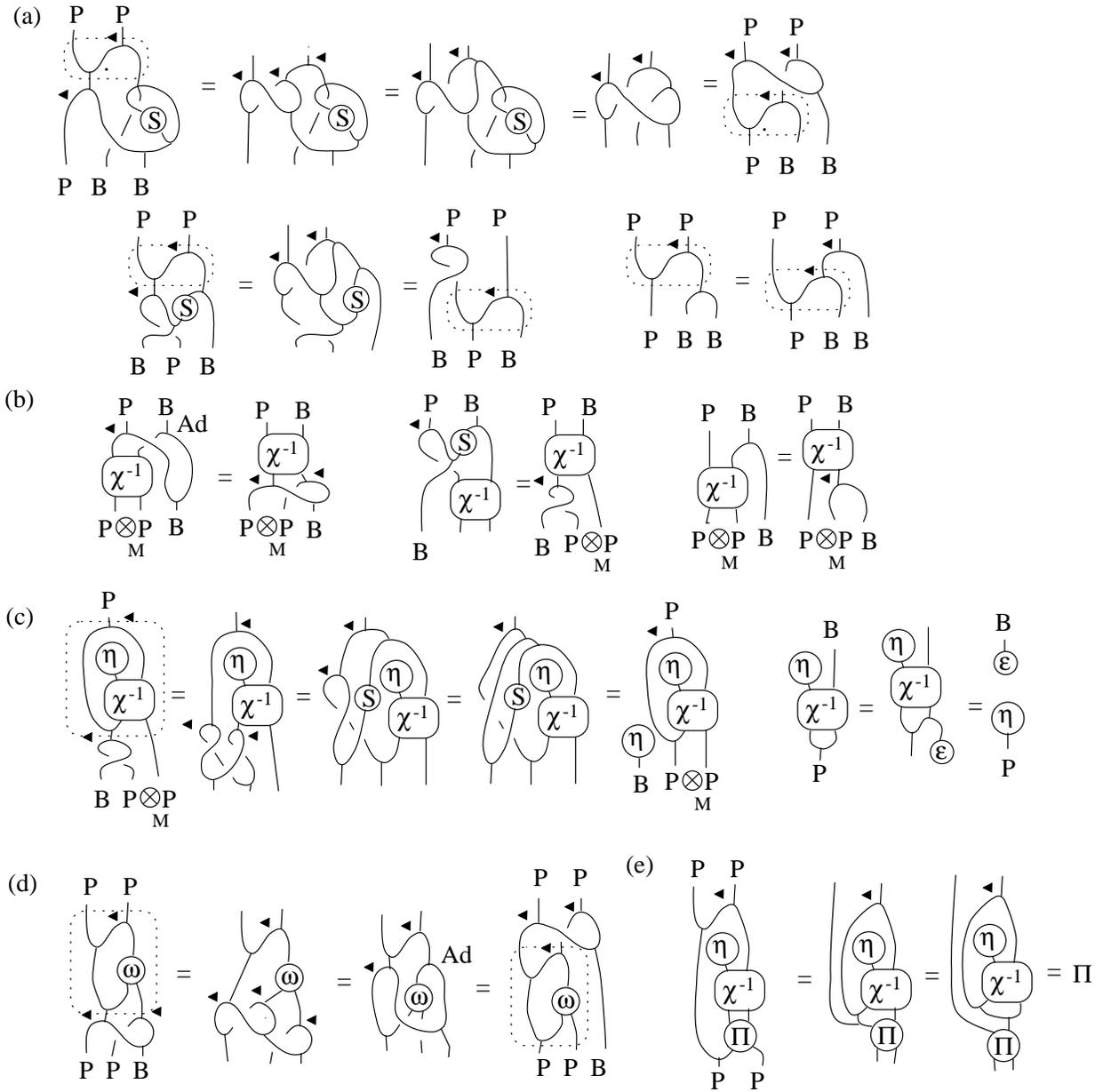}\]
\caption{Correspondence of projections and connections:  (a) covariance properties of $\chi$ and (b)--(c) of $\chi^{-1}$ needed for construction (d) of $\Pi$  from $\omega$ and (e) its reconstruction}
\end{figure}
\proof This is shown in Figure~3. Part the result (not the explicit reconstruction of $\Pi$ in the 1-1 correspondence) is in \cite{BrzMa:coa} in an algebraic form as an example of coalgebra gauge theory. Parts (a)-(c) show first some covariance properties of $\tilde\chi$ and $\chi^{-1}$. Part (a) shows that $\tilde\chi:P\tens P\to P\tens B$ (box on right) is an intertwiner for the tensor product coactions of $B$, where $B$ acts on $B$ by the braided adjoint coaction (box on left) from \cite{Ma:lin}. We use the homomorphism property of $\cora$, the comodule property of $\cora$, the antipode axioms to cancel the loop with $S$, and finally the comodule property in reverse. For other properties if $\tilde\chi$ shown are immediate in a similar way. We deduce without any work the corresponding properties for $\chi^{-1}$, shown in part (b). Part (c) deduces some further properties of the combination involving $\chi^{-1}\circ(\eta\tens\id)$, using the comodule homomorphism property, part (b), the comodule property and the antipode axioms. The first result is that the combination (in box) lies in $M\tens_MP$. The second  result in part (c) follows by writing the product in $P$ as an application of $\chi$. Part (d) proves that $\Pi$ defined from $\omega$ is indeed an intertwiner $\Omega^1P\to\Omega^1P$ for the coaction of $B$. We use the homomorphism property of $\cora$, (C2) and part (a). The other list properties of $\Pi$ are immediate from its form as a composite of $\omega$ and $\tilde\chi$, given (C1). Finally, given an idempotent $\Pi$ with these properties, we define $\omega$ as stated in the proposition. This is well defined because $\Pi$ vanishes on $P(\Omega^1M)P$ and hence, in particular, on $P(dM)P$. Here the restriction to $\ker\eps$ of $\chi^{-1}(\eta\tens \id)$ factors through (the projection to $P\tens_MP$ of) $\Omega^1P$ because of the second result in part (c). Then (C1) holds in view of the assumption $\tilde\chi\circ\Pi=\tilde\chi$, and (C2) in view of the assumed covariance of $\Pi$ under $B$ and part (b). That these constructions are inverse in one direction is trivial from their form. The proof in the other direction, defining $\omega$ from $\Pi$ and then computing its associated projection, is shown in part (e). We use that $\Pi$ is assumed a left $P$-module morphism, and then part (c). The reconstruction of $\Pi$ in this way is new even in the quantum group case, being covered somewhat implicitly in \cite{BrzMa:gau}. \endproof

The condition $\tilde\chi\circ\Pi=\tilde\chi$ can also be cast as $\ker\Pi=P(\Omega^1M)P$ (given that $\Pi$ is already assumed to be zero on this), as in the classical setting. So a projection provides an equivariant complement to the horizontal forms. Also, let $\omega$ be a connection on $P$ and $\alpha:B\to \Omega^1P$ an intertwiner as in (C2) such that $\tilde\chi\circ\alpha=0$ and $\alpha\circ\eta=0$. Then $\omega+\alpha$ is also a connection
on $P$, and  the difference of any two connections is of this form. Hence we identify $\CA(P,B)$, the space of connections on $P$, as an affine space.  Finally, as discovered in \cite{BrzMa:gau}, not all connections come locally from the base when our algebras are noncommutative. This is due to the distinction between $P(\Omega^1M)P$ and $(\Omega^1M)P$ in the noncommutative case;
 
(C3) A connection is said to be {\em strong} if $(\id-\Pi)\circ d:P\to \Omega^1P$ factors through $(\Omega^1M)P$. 

This is a braided version of the condition recently developed by P.~Hajac in \cite{Haj:str} for quantum group gauge theory. We will use it especially  Section~3. In terms of $\omega$ the condition is
\[ \epsfbox{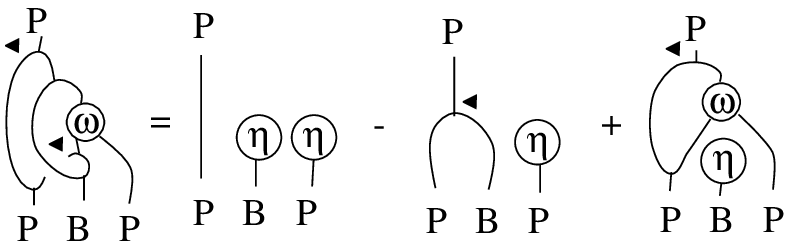}\]
in the case of the universal differential calculus. The collection $\CA_s(P,B)$ of strong connections also forms an affine space. Differences $\alpha$ between strong connections obey the above diagram with $\omega$ replaced by $\alpha$ and without the terms not involving $\alpha$.

Finally, we consider the natural gauge equivalence of these various data. 
Let $P,B$ be a braided principal bundle. We define  a global gauge transformation to be a morphism $\Gamma:B\to P$ such that:

(G1)  $\Gamma$ is convolution-invertible and  obeys $\Gamma\circ\eta=\eta_P$.

(G2)  $\cora\circ\Gamma=(\Gamma\tens\id)\circ\circ\Ad$, where $\cora$ is the given coaction on $P$ and $\Ad$ the braided  adjoint coaction from \cite{Ma:lin}.

\begin{propos} Gauge transformations $\Gamma$  are in 1-1 correspondence with invertible morphisms $\Theta:P\to P$ which are left $M$-module morphisms and  intertwiners for the coaction of $B$ obeying $\Theta\circ\eta_P=\eta_P$, via
\[  \epsfbox{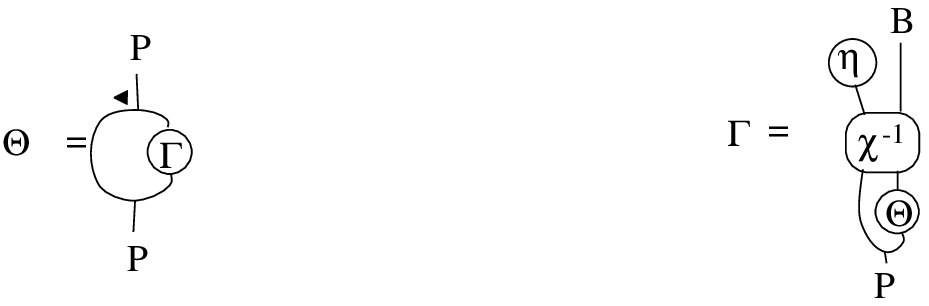}\]
Moreover, the collection $\CG(P,B)$ of gauge transformations  forms a group under convolution product, and $\Theta_{\Gamma*\Gamma'}=\Theta_{\Gamma'}\circ\Theta_{\Gamma}$ represents it on $P$. 
\end{propos}
\begin{figure}
\[ \epsfbox{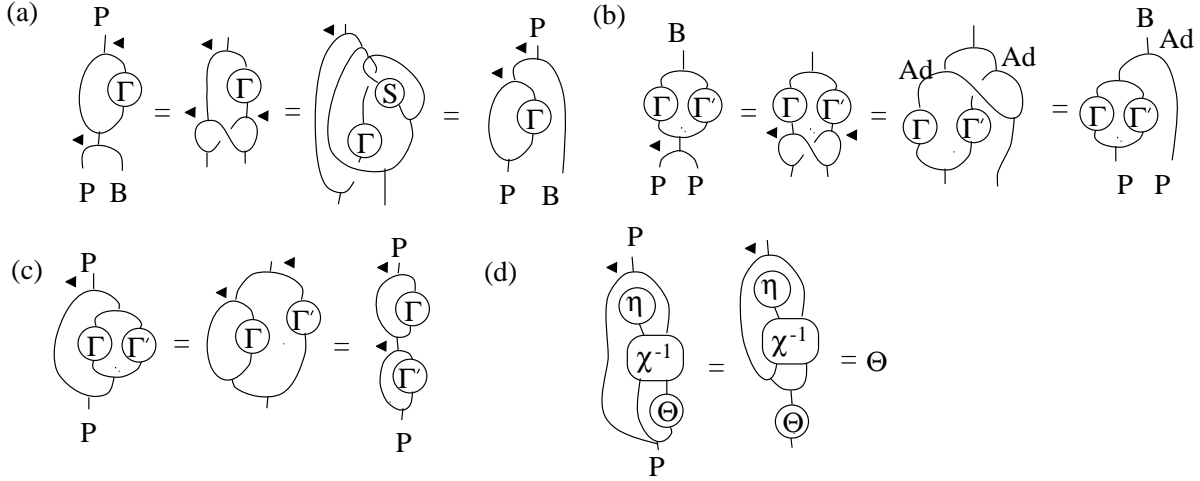}\]
\caption{Construction (a) of bundle morphism $\Theta$ from gauge transformation $\Gamma$ and (b) its reconstruction}
\end{figure}
\proof This is shown in Figure~4. This is actually very similar to (a simpler version of)  the proof for $\Pi$ and $\omega$ in the preceding proposition. Part (a) checks that $\Theta$ is an intertwiner for the coaction of $B$. We use the homomorphism property, the comodule property and (G2), cancel the antipode loop and then the comodule property in reverse. The $M$-module and unity properties of $\Theta$ are immediate. Part (b) checks closure under the convolution product and follows at once from the fact (used in the final equality) that the braided adjoint coaction respects the braided coproduct\cite{Ma:lin}. Hence the convolution $\Gamma*\Gamma'$ also obeys (G2). Part (c) shows that the $\Theta$ construction represents the convolution product by composition. We use the comodule property followed by  part (a). Conversely, starting from $\Theta$  it is immediate from the properties if $\chi^{-1}$ in Figure~3(b) that $\Gamma$ obeys (G2). That these constructions are inverse is also immediate in one direction from their form. In the other direction, starting from a morphism $\Theta:P\to P$, we define $\Gamma$ as stated in the proposition and reconstruct $\Theta$. This is shown part (d). We use Figure~3(c) and that $\Theta$ is assumed to be a left $M$-module morphism. \endproof

Such global gauge transformations have been considered previously only in the quantum group case, in \cite{Brz:tra}. Note, however, that $\Theta$ is not a bundle automorphism in a natural sense because it need not respect the algebra structure of $P$. Rather, we think of it as a {\em bundle transformation}  $P\to P^\Gamma$ where $P^\Gamma$ has a new product $\cdot_\Gamma=\Theta\circ\cdot\circ(\Theta^{-1}\tens\Theta^{-1})$ and forms a comodule-algebra under the same coaction $\cora$. We say that $P^\Gamma, B$ and $P,B$ are {\em globally gauge equivalent}.

That $\CG$ modifies the algebra structure of $P$ (isomorphically) is an interesting complication arising from its non-commutativity. Apart from this, it acts as well on connections (preserving the strong connections) by
\[ \omega^\Gamma=(\Theta\tens\Theta)\circ\omega.\]
This is arranged so that when we compute $\Pi^\Gamma:\Omega^1P^\Gamma\to \Omega^1P^\Gamma$ using $\omega^\Gamma$ and $\cdot_\Gamma$ (in $\chi$ and the definition of $\Omega^1P^\Gamma$), we have the commuting square $(\Theta \tens\Theta )\circ\Pi=\Pi^\Gamma\circ(\Theta \tens\Theta )$. This means that $(P^\Gamma,\omega^\Gamma,B,\cora)$ and $(P,\omega,B,\cora)$ are the same abstract connection and bundle after allowing for the algebra isomorphism $\Theta :P\isom P^\Gamma$. This point of view of $\Theta$ as a bundle transformation appears to be new even in the quantum group case.

\section{Trivial Braided Principal Bundles and Gauge Fields}

In this section we study a class of examples of braided principal bundles provided by the following data:

(T1)  $P$ a braided right $B$-comodule algebra, with $M=P^B$, i.e. (P1) as above.

(T2) A morphism $\Phi:B\to P$ such that $\cora\circ\Phi=(\Phi\tens\id)\circ\Delta$ and $\Phi\circ\eta=\eta_P$.

(T3) A morphism $\Phi^{-1}:B\to P$ such that $\cdot_P\circ(\Phi\tens\Phi^{-1})\circ\Delta=\cdot_P\circ(\Phi^{-1}\tens\Phi)\circ\Delta=\eta_P\circ\eps$. This morphism, if it exists, is uniquely determined by $\Phi$ and is called its {\em convolution inverse}.

\begin{propos} Given the data $P,B,\Phi$ obeying (T1)-(T3) we have 
$M\tens B\isom P$ as objects, and a braided principal bundle structure via the morphisms 
\[\epsfbox{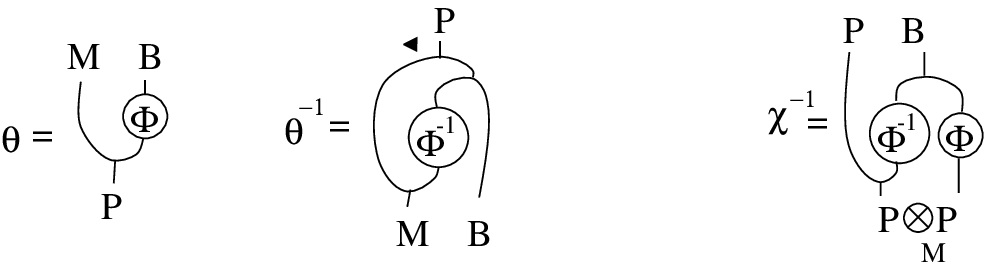} \]
We call this the {\em trivial principal bundle associated to $\cora,\Phi$}.
\end{propos}
\begin{figure}
 \[  \epsfbox{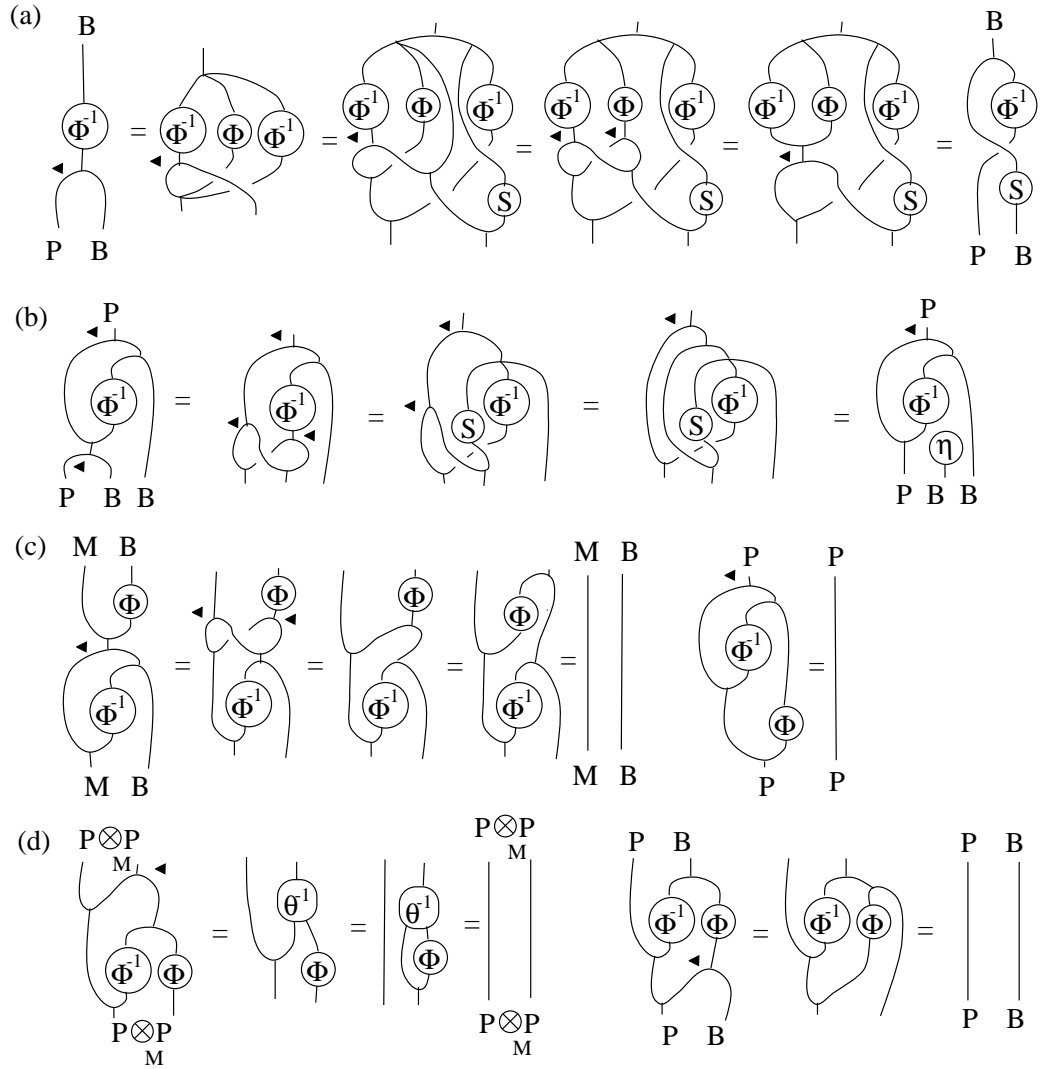}\]
\caption{Construction of trivial principal bundles showing (a) covariance of $\Phi^{-1}$ (b) proof that $\theta^{-1}$ factors through  $M\tens B$ (c) proof that $\theta,\theta^{-1}$ are inverse (d) proof that $\chi,\chi^{-1}$ are inverse}
\end{figure}
\proof This is in \cite{BrzMa:coa} in an algebraic form, as an example of coalgebra gauge theory. We complement this now with the braid-diagrammatic proof in Figure~5. Part (a)  begins with the covariance property of $\Phi^{-1}$. We insert a trivial $\Phi,\Phi^{-1}$ loop, and then a trivial antipode loop. We  then use the covariance of $\Phi$, and the  homomorphism property of the coaction. We then cancel the $\Phi^{-1},\Phi$ loop. Part (b) verifies that $\theta^{-1}$ factors through $M\tens B$ as required. We use the homomorphism property, part (a) and the comodule property. We then cancel and antipode loop. Part (c) verifies that $\theta,\theta^{-1}$ are indeed inverse. We use the homomorphism property, that the action is trivial on $M$ and the covariance of $\Phi$. The inverse from the other side it immediate by cancelling the $\Phi^{-1},\Phi$ loop. Part (d) verifies that $\chi,\chi^{-1}$ are inverse. We use part (b), allowing to move the product in $P\tens_MP$, and part (c). The inverse from the other side is immediate from the covariance of $\Phi$. The quantum group case is in \cite{BrzMa:gau}. \endproof

\begin{example} Let $M$ be an algebra and $B$ a braided group in our braided category. Then the {\em braided tensor product algebra} $M\und\tens B$ is a trivial quantum principal bundle with 
\[ \epsfbox{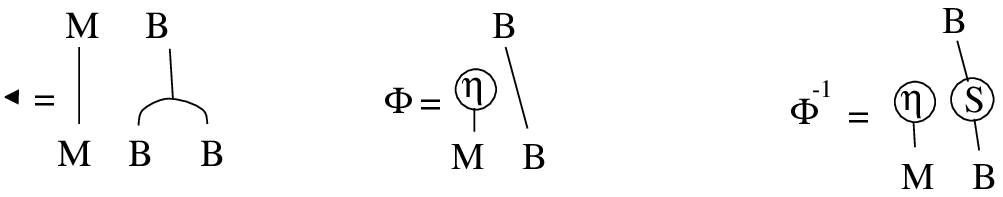}\]
\end{example}
\begin{figure}
 \[  \epsfbox{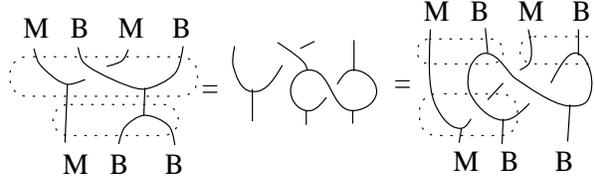}\]
 \caption{Proof of braided comodule algebra structure needed for tensor product principal bundle}
\end{figure}

\proof This is shown in Figure~6. The coaction $\cora$ shown in the lower box on the left is the (braided) tensor product coaction of the trivial coaction on $M$ and the right regular coaction on $B$. We check that the braided tensor product algebra structure $M\und\tens B$ as shown in the upper box on the left forms a braided comodule algebra under it, using the homomorphism property of the coproduct. The other facts needed to obey (T1)-(T3) are obvious. \endproof

The  main result which justifies the notion of trivial braided principal bundles is that associated to every trivialisation is a class of connections of the form 
\[ \omega_{A,P,\Phi}=\Phi^{-1}*d\Phi+\Phi^{-1}*A*\Phi,\] 
where  
 $A:B\to \Omega^1M$ is any morphism such that $A\circ\eta=0$. We start with an abstract characterisation of the class of connections which arise this way.

\begin{propos} Let $P,B,\Phi$ be a trivial braided principal bundle.  Then strong connections $\omega$  are in 1-1 correspondence with morphisms $A:B\to \Omega^1M$ such that $A\circ\eta=0$, via
\[ \epsfbox{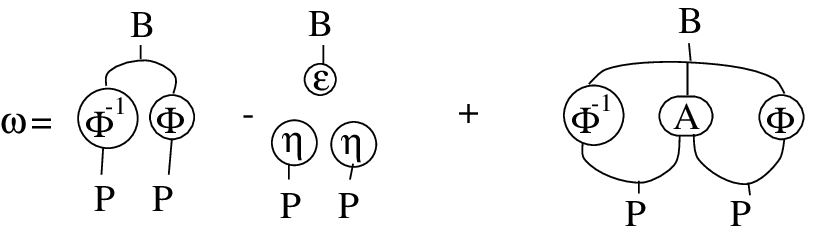}\]
We call such $A$ {\em braided gauge fields or local connections}.
\end{propos}
\begin{figure}
 \[  \epsfbox{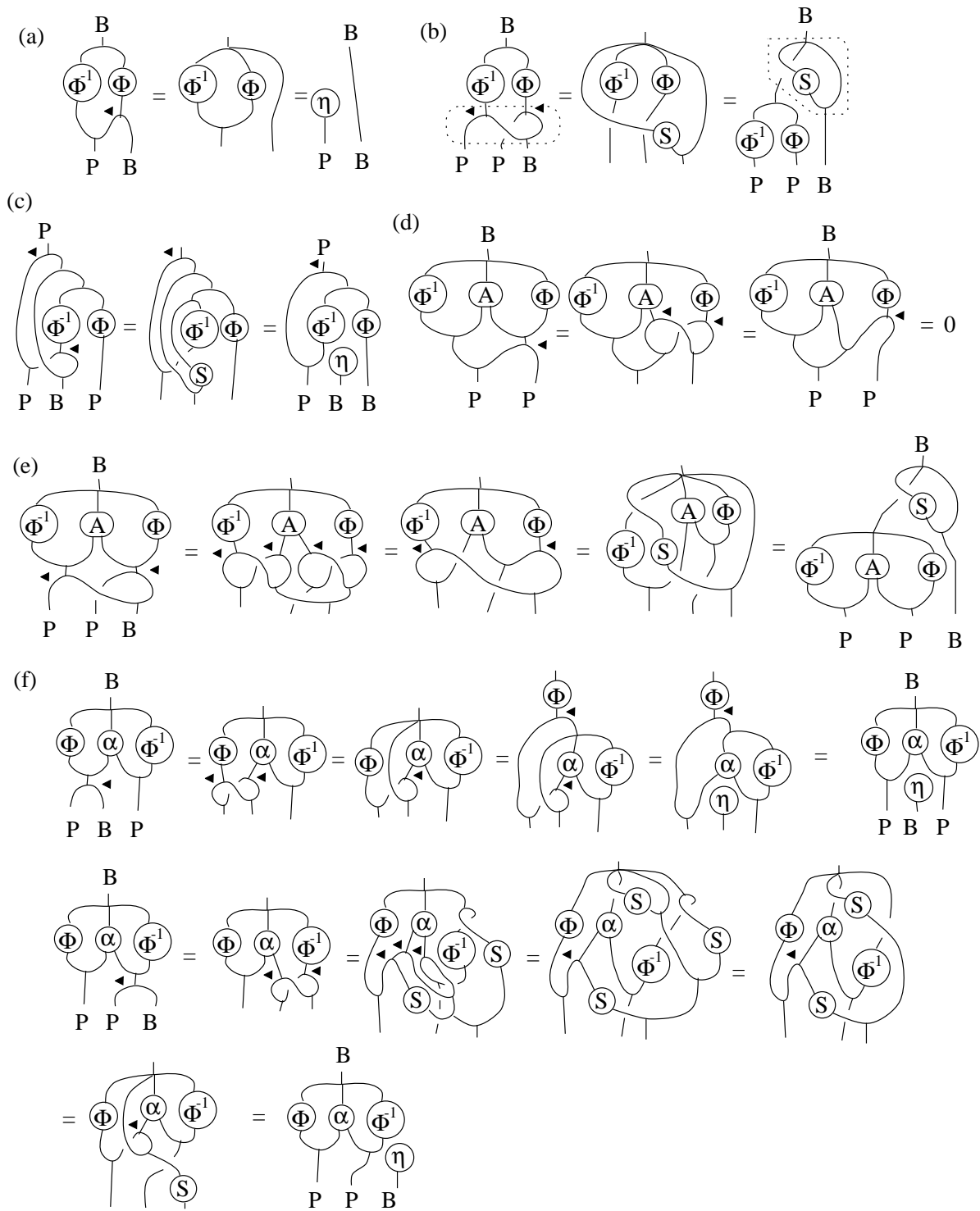}\]
\caption{Construction of connections on trivial braided principal bundles; (a)--(c) for trivial connection, (d)--(e) additional part of connection defined by braided gauge field $A$. Conversely (f) from a strong connection we construct a gauge field}
\end{figure}
\proof The proof is in Figure~7. Part of the result (not the strongness condition and not the  reconstruction of $A$ from $\omega$) is  in \cite{BrzMa:coa} in  an algebraic form as an example of coalgebra gauge theory. Part (a) applies $\tilde\chi$ to the main part of the trivial connection $\omega_{0,P,\Phi}=\Phi^{-1}*d\Phi$. We use the covariance of $\Phi$. Hence (C1) follows. Part (b) uses the covariance of $\Phi^{-1},\Phi$ to obtain (C2) for the trivial connection. Part (c) verifies (C3) for the trivial connection. We use the covariance of $\Phi^{-1}$ from Figure~5(a) and then cancel an antipode loop. Part (d) applies $\tilde\chi$ to the part of the connection coming from a gauge field $A:B\to \Omega^1M$. We use that $\cora$ is trivial on $M$. Hence (C1) holds. Part (e) uses the homomorphism property of the coaction, the triviality on $M$ and the covariance of $\Phi^{-1},\Phi$ to obtain (C2). The final steps are just the same as for the trivial connection. Similarly, (C3) holds by just the same steps as for the trivial connection. Finally, we consider the difference $\alpha$ between a given connection and the trivial one, and define $A=\Phi*\alpha*\Phi^{-1}$. Part (f) verifies that this indeed lies in $\Omega^1M$ provided our connection obeys (C3). We use the comodule homomorphism property, the covariance of $\Phi$, and (after coassociativity) the covariance of $\Phi$ in reverse. Finally, we apply (C3) in the form for $\alpha$ and the covariance of $\Phi$ again. This shows that the left hand output of $A$ factors through $M$. For the right hand output the proof is more complicated. We use the homomorphism property, the covariance of $\Phi^{-1}$ from Figure~5(a) and also insert a trivial antipode loop. We then use the intertwiner property (C2) and cancel the resulting antipode loop. Finally, we use the braided antihomomorphism property of the braided antipode\cite{Ma:tra} to obtain an expression similar to the third diagram in the proof for the left hand output. The rest of the proof follows that case. The construction of $\omega$ from $A$ in the quantum group case is in  \cite{BrzMa:gau} and the converse in this case is in \cite{Haj:str}. \endproof

Given a trivial braided principal bundle $P,B,\Phi$ we define its group of {\em local gauge transformations} as the convolution-invertible morphisms $\gamma:B\to M$ such that $\gamma\circ\eta=\eta_M$. These form an ordinary group  under convolution and act transitively on the collection of trivialisations by  $\Phi^\gamma\equiv \gamma*\Phi=\cdot_P\circ(\gamma\tens\Phi)\circ\Delta$, i.e. this is also a trivialisation and any two trivialisations of the same bundle are related by such a $\gamma$ (this follows at once from Proposition~4.1 in the next section). We say that the corresponding trivial principal bundles are {\em locally gauge equivalent}. If $A$ is a gauge field, we define its gauge transform 
\[ A^\gamma=\gamma^{-1}*A*\gamma+\gamma^{-1}*d\gamma\]
which is such that $\omega_{A^\gamma,P,\Phi}=\omega_{A,P,\Phi^\gamma}$ is the same connection on $P$ when defined via the gauge transformed trivialisation. This is the passive view of gauge transformations. These steps work in just the same way as in the quantum group case in \cite{BrzMa:gau}.

\begin{propos} Let $P,B,\Phi$ be a trivial braided principal bundle. Then global gauge transformations $\Gamma$ are in 1-1 correspondence with local gauge transformations $\gamma$, via
\[ \epsfbox{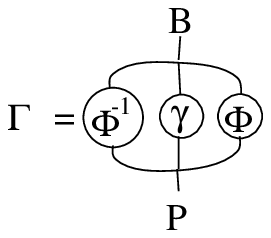}\]
The corresponding gauge transformed trivialisation is $\Phi^\Gamma=\Theta\circ\Phi=\gamma*\Phi$. Moreover, $(\omega_{A,P,\Phi})^\Gamma=\omega_{A,P^\Gamma,\Phi^\Gamma}$.
\end{propos}
\begin{figure}
\[ \epsfbox{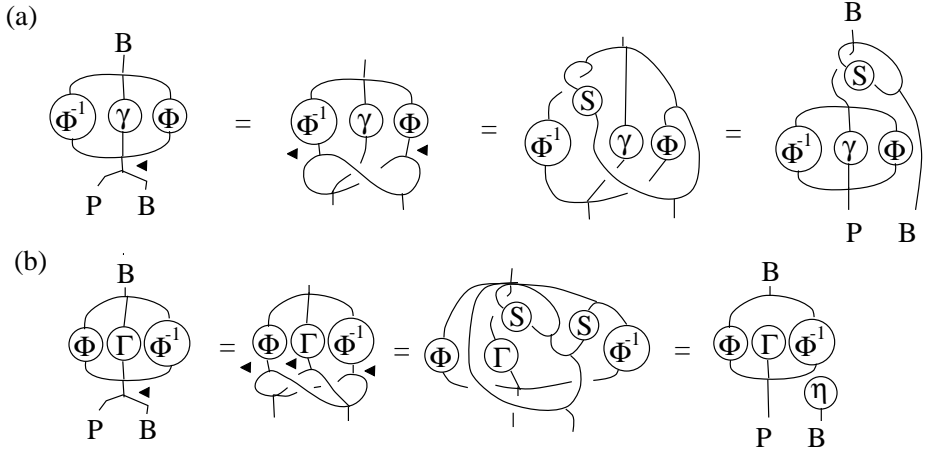}\]
\caption{Equivalence of local and global gauge transformations on a trivial bundle}
\end{figure}
\proof This is shown in Figure~8. The proof is similar to (a simpler version) of the proof for gauge fields in the preceding proposition. Part (a) verifies that $\Gamma$ defined from $\gamma$ obeys (G2). We use the comodule homomorphism property (with trivial coaction on $M$), and covariance of $\Phi^{-1},\Phi$. Conversely, part (b) verifies that $\gamma=\Phi*\Gamma*\Phi^{-1}$ factors through $M$. We use the comodule homomorphism property and covariance of $\Phi,\Gamma,\Phi^{-1}$ and cancel the resulting antipode loops. The remaining assertions are equally straightforward. The quantum group case is a version of identical arguments in \cite{BrzMa:gau} for gauge fields and pseudotensorial forms \endproof

Finally, we note that our isomorphism $\theta$ in Proposition~3.1 puts every
trivial braided principal bundle with chosen trivialisation into a canonical form built on $M\tens B$. As such, trivial braided principal bundles correspond to braided convolution-invertible cocycle cross products. 

\begin{figure}
\[ \epsfbox{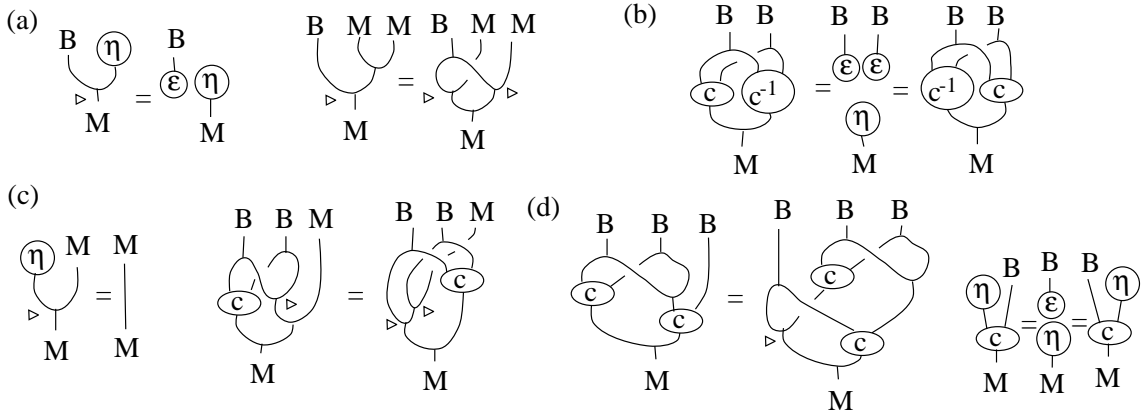}\]
\caption{Braided 2-cocycle data}
\end{figure}
\begin{propos} Let $c:B\tens B\to M$ and $\la:B\tens M\to M$ form a braided 2-cocycle in the sense shown in Figure~9. Then the braided cocycle cross product  $P=M{}_c\lcross B$ built on $M\tens B$ with product, coaction and trivialisation
\[ \epsfbox{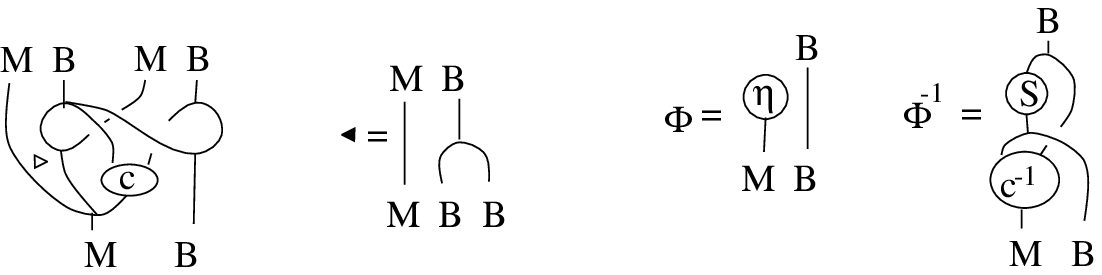}\]
forms a trivial braided principal bundle. Conversely, every trivial braided principal bundle with chosen trivialisation is isomorphic via $\theta$ to a braided cocycle cross product bundle, with
\[ \epsfbox{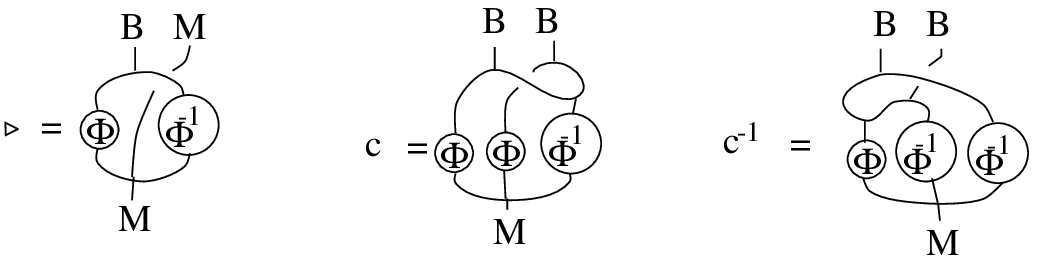}\]
\end{propos}
\begin{figure}
\[ \epsfbox{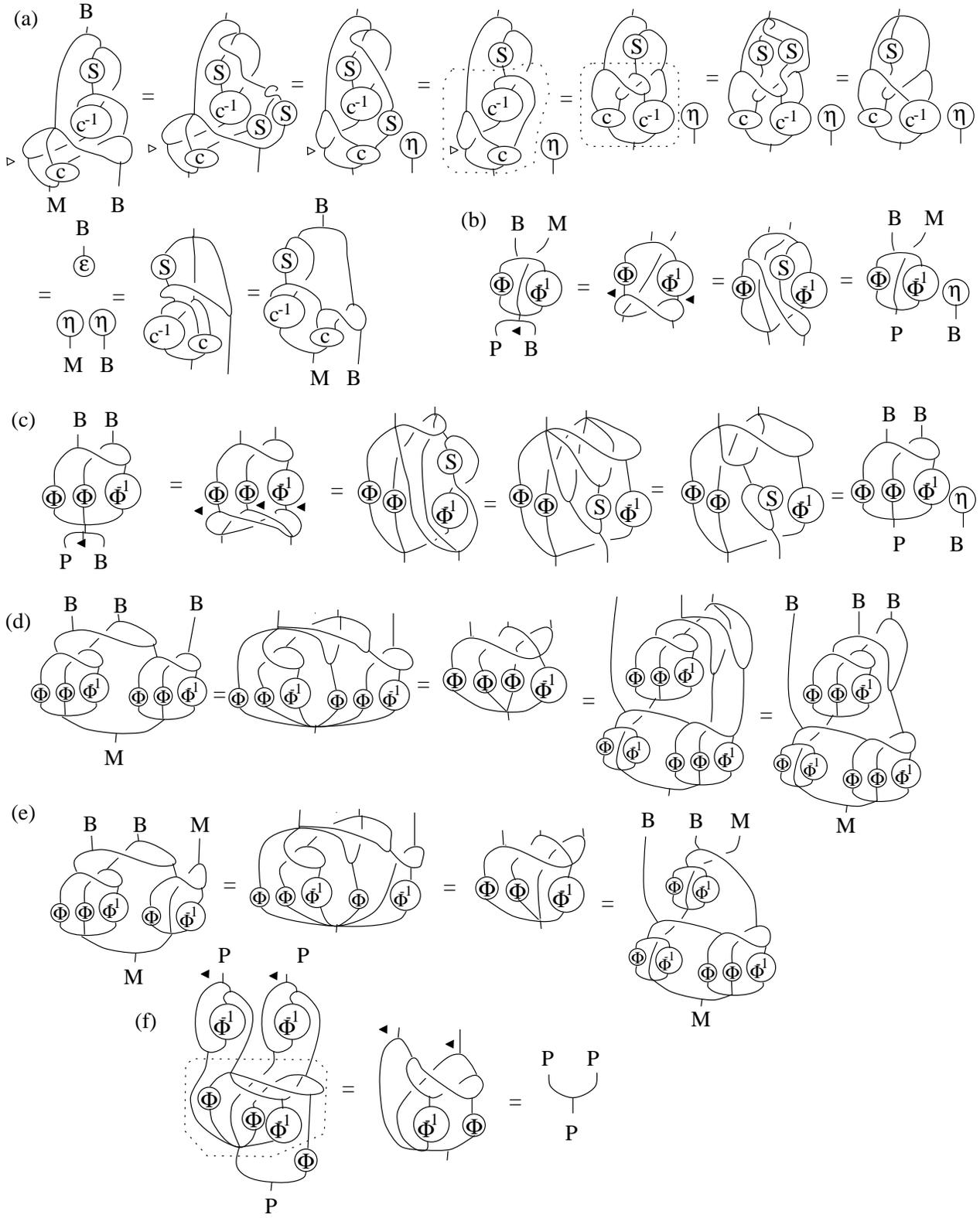}\]
\caption{(a) Proof of $\Phi^{-1}$ for cocycle cross product bundle, (b)--(e) construction of cocycle data from a bundle and (f) isomorphism with $P$}
\end{figure}
\proof Part of the result (not the construction of a bundle from a cocycle) is in \cite{Brz:coa} in an algebraic form as an example of cross products in coalgebra gauge theory; we provide direct  braid-diagrammatic proofs. In fact, the proof that product stated on $M\tens B$ (defined by cocycle $c,\la$) is associative  follows just the same lines as given in detail for cross products by braided groups (without
cocycle) in \cite{Ma:bos}. That the stated coaction makes $M_c\lcross B$ a braided comodule algebra follows the same argument as in Figure~6. The proof that $\Phi,\Phi^{-1}$ as stated are inverse is shown in Figure~10(a). On the left we show $\Phi*\Phi^{-1}$ computed with the product of $M_c\lcross B$. We then use (twice) the braided-antimultiplicativity of the braided antipode and cancel a resulting antipode loop. We then use the braided-antimultiplicativity property in reverse. Next we use the cocycle axiom from Figure~9(d) but massaged in the form shown in the two boxes. Equality of the two boxes is equivalent to Figure~9(d) after convolving twice with $c$. The sense in which $c,c^{-1}$ are inverse is in Figure~9(b), namely under the convolution product on morphisms $B\und\tens B\to M$, where $B\und\tens B$ is the braided tensor coalgebra. After this, we use the braided antimultiplicativity of $S$ one more time, cancel an antipode loop and $c,c^{-1}$. The computation for $\Phi^{-1}*\Phi$ is more immediate. We consider now the converse direction, starting from a trivial braided principal bundle $P,B,\Phi$. Part (b) verifies that $\la$ as stated factors through $M$. We use the comodule homomorphism property, the covariance of $\Phi,\Phi^{-1}$ and cancel the resulting antipode loop. Part (c) similarly verifies that $c$ as stated factors through $M$. After covariance of $\Phi,\Phi^{-1}$ we use the homomorphism property of the coproduct, braided antimultiplicativity of $S$ and the coproduct homomorphism property in reverse,  and can cancel the resulting antipode loop. Part (d) verifies the cocycle axiom in Figure~9(d). We use the coproduct homomorphism property twice, allowing us to cancel $\Phi^{-1},\Phi$. We can then insert several cancelling $\Phi^{-1},\Phi$ loops and use the coproduct homomorphism property again. The proof for the cocycle axiom in Figure~9(c) is very similar and shown in part (e). The remaining cocycle requirements are more immediate. Finally, part (f) shows that the braided cross product with these cocycles is isomorphic under $\theta$ from Proposition~3.1 to $P$. The product for this cocycle simplifies to the expression in the box on cancelling $\Phi,\Phi^{-1}$. Further cancellation recovers the product of $P$. \endproof

 As we change our choice of trivialisation, we will change the corresponding cocycle. We also have an active point in which the trivialisation is fixed and the bundle itself changes to $P^\gamma$ under the bundle automorphism $\Theta$ induced by local gauge transformation $\gamma$. 
    
\begin{propos} Let $\gamma:B\to M$ be a local gauge transformation. If a trivial bundle $P,B,\Phi$ has canonical form $M{}_c\lcross B$ then $P,B,\Phi^\gamma$ and $P^{\gamma^{-1}},B,\Phi$ have the canonical form  $M{}_{c^\gamma}\lcross B$, where \[ \epsfbox{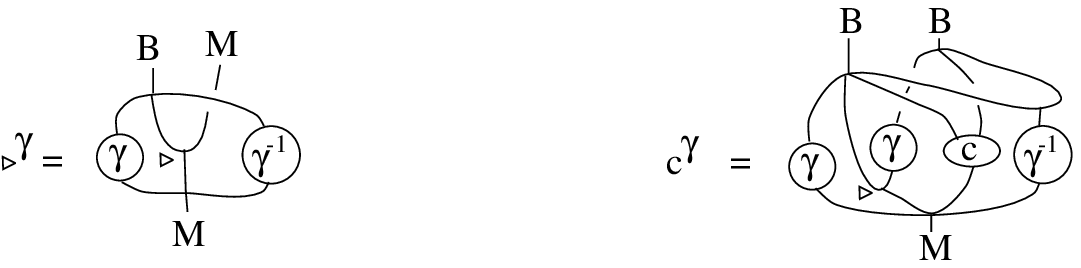}\]
form another braided 2-cocycle, the {\em gauge transform} of $c,\la$. 
In particular, $(M{}_c\lcross B)^{\gamma^{-1}}=M{}_{c^{\gamma}}\lcross B$ and the induced $\Theta_{\gamma^{-1}}:M_c\lcross B\to M_{c^\gamma}\lcross B$ is an isomorphism of braided comodule algebras. 
\end{propos}
\begin{figure}
\[  \epsfbox{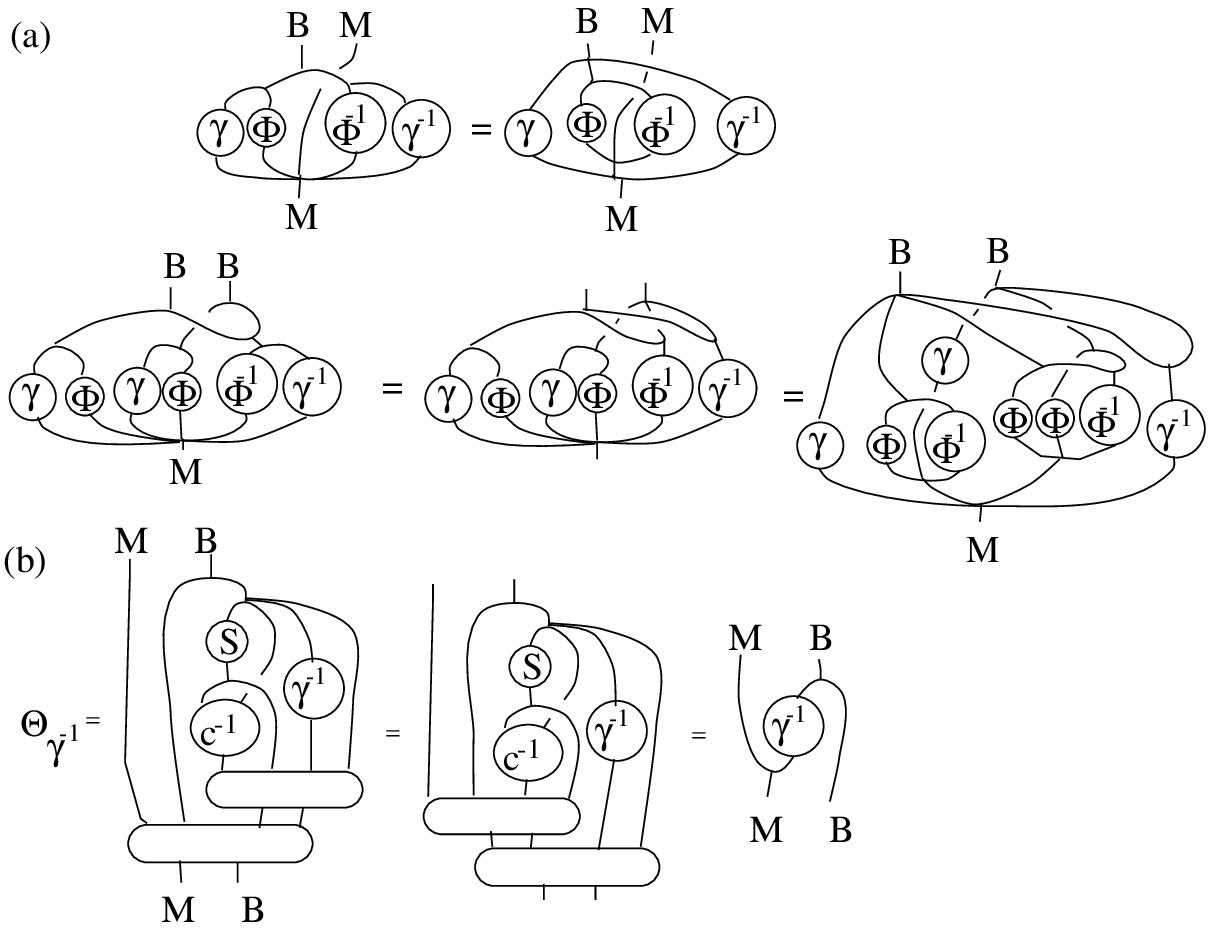}\]
\caption{(a) Gauge transformation of 2-cocycle induced by change of trivialisation and (b) induced bundle morphism}
\end{figure}
\proof  The proof is in Figure~11. Part (a) shows that the cocycle  computed  with $\Phi^\gamma=\gamma*\Phi$  gives $c^\gamma,\la^\gamma$, where $c,\la$ are the cocycles computed with respect to $\Phi$. We use the coproduct homomorphism property and insert a trivial $\Phi^{-1},\Phi$ loop for the $c^\gamma$ part. Since $\Theta\circ\Phi=\gamma*\Phi$ from Proposition~3.4 it is immediate that the cocycle computed from $P^{\gamma^{-1}},\Phi$ (using the product of $P^{\gamma^{-1}}$)  is also   $c^\gamma,\la^\gamma$.  In part (b) we compute the associated  bundle isomorphism $\Theta_{\gamma^{-1}}:P\to P^{\gamma^{-1}}$ in the particular case  $M_c\lcross B\to M_{c^\gamma}\lcross B$. We use the definitions in Propositions~2.2 and~3.4, and $\Phi,\Phi^{-1}$ for the cross product bundle from Proposition~3.5. The boxes are the product in $M{}_c\lcross B$. We then recognise a product $\Phi*\Phi^{-1}$ and cancel it. One may easily verify directly that it is an isomorphism of braided comodule algebras fixing $M$. \endproof

We see that equivalence classes of trivial bundles correspond to braided  2-cocycles `up to coboundary' in the sense of up to gauge transform of the cocycles, i.e. to braided non-Abelian 2-cohomology.  This is just as in the quantum group case in \cite{Doi:equ}\cite{Ma:clau}.
Note also that this cohomology has no analogue in classical geometry, where commutativity of our algebras forces the cocycle to be trivial.   
 
\begin{example} In the case of a braided tensor product bundle the strong connections and global gauge transformations take the form 
\[ \epsfbox{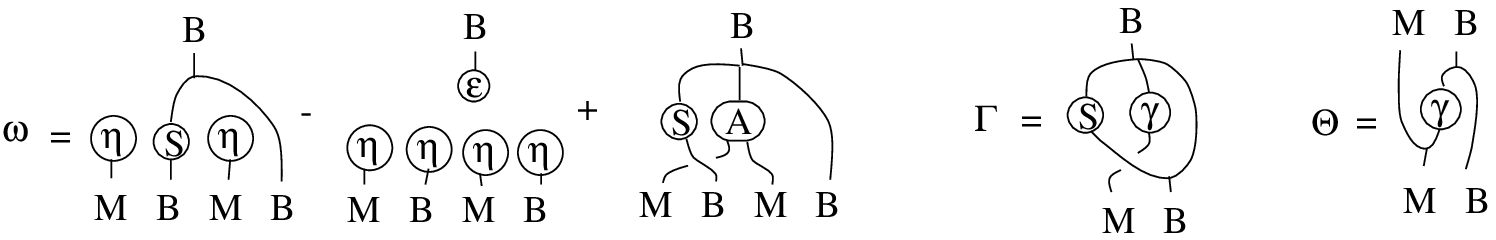}\]
The gauge transformed bundle $(M\und\tens B)^{\gamma^{-1}}$ via $\Theta_{\gamma^{-1}}$ is no longer a braided tensor product but has the form $M{}_c\lcross B$ with product and associated cocycle (a coboundary)  
\[ \epsfbox{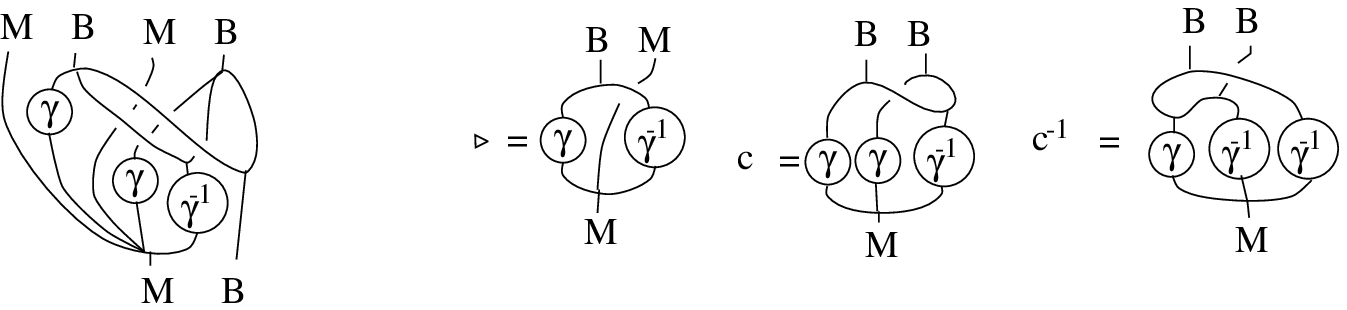}\]
\end{example}
\proof We use the particular form of the braided tensor product algebra $M\und \tens B$ to obtain the form of $\omega$ and $\Gamma$ and corresponding bundle morphism $\Theta$. For the product of $(M\und\tens B)^{\gamma^{-1}}$ we easily compute $\Theta_{\gamma^{-1}}\circ\cdot\circ(\Theta_\gamma\tens\Theta_\gamma)$ using form of $\Theta$ stated. By the preceding proposition, it must have the cocycle form  shown. \endproof 

\section{Covariant Derivative and Associated Braided Fiber Bundles}

Let $V$ be a right $B$-comodule with coaction $\beta$. We consider fields with values in `the braided space with coordinate ring' $V$, i.e. mapping from $V$. For a full geometrical picture $V$ should be a braided comodule algebra but most results in this section do not actually need this. As for the quantum group case in \cite{BrzMa:gau}, we work globally on $P$ and define a {\em braided pseudotensorial $n$-form} on a braided principal bundle $P,B$ as a morphism $\Sigma: V\to \Omega^nP$ which is equivariant, i.e. intertwines the coactions on $V,P$.  Examples of pseudotensorial forms already encountered above are  trivialisations $\Phi$, which are pseudotensorial forms $\Phi:B_R\to P$, global gauge transformations $\Gamma:B_{\Ad}\to P$ and differences of connections $\alpha:B_{\Ad}\to \Omega^1P$. Here $B_R$ denotes $B$ under the right regular coaction $\Delta$ and $B_\Ad$ denotes it with the braided adjoint coaction from \cite{Ma:lin}. The former is always a braided comodule algebra, while the latter is\cite{Ma:lie} a braided comodule algebra whenever $B$ is braided-commutative with respect to $\Ad$ in the sense below.

In fact, the forms that we are interested in are not these pseudotensorial ones but forms coming from the base $M$. The approach in \cite{BrzMa:gau} is to attempt to restrict to `strongly tensorial' ones, namely $\Sigma$ which factor through $(\Omega^nM)P$. The general picture is not well understood, but this is indeed what happens at least for trivial braided principal bundles. 

\begin{propos} Let $P,B$ be a trivial braided principal bundle. Then strongly tensorial forms $\Sigma$ are  in 1-1 correspondence with morphisms $\sigma:V\to \Omega^nM$, via
\[ \epsfbox{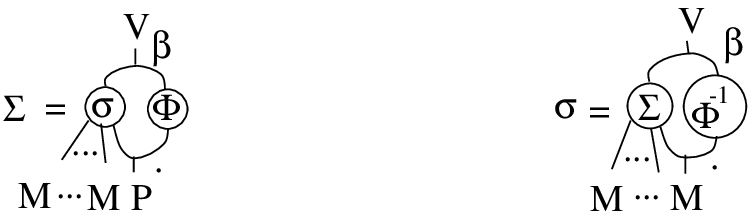}\]
We write $\Sigma=\sigma*\Phi$ and $\sigma=\Sigma*\Phi^{-1}$ as an extension of the convolution notation. Morphisms $\sigma$ are called {\em local sections}.
\end{propos}
\begin{figure}
\[ \epsfbox{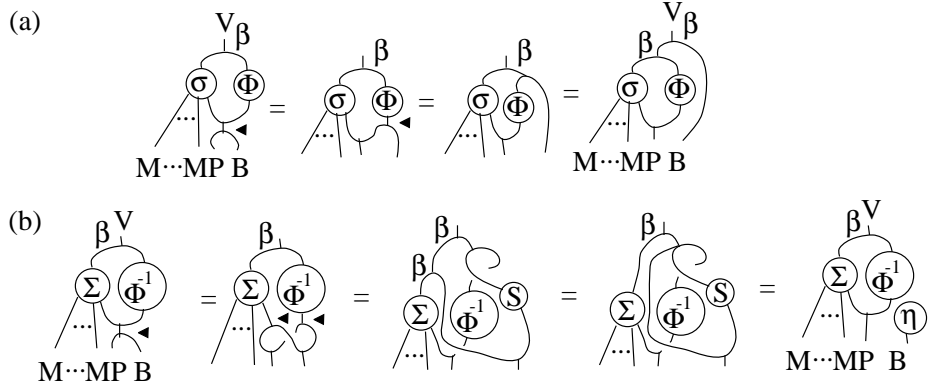}\]
\caption{Equivalence of strongly tensorial forms and local sections on a trivial bundle}
\end{figure}
\proof This is shown in Figure~12. Part (a) shows that $\Sigma$ defined from any morphism $\sigma:V\to M$ is pseudotensorial. We use the comodule homomorphism property, covariance of $\Phi$ and the comodule property. It is manifestly strongly tensorial. Part (b) shows that $\sigma$ defined from a pseudotensorial form $\Sigma$ factors through $\Omega^nM$. We use the comodule homomorphism property, covariance of $\Phi^{-1}$ and the comodule property, and cancel the resulting antipode loop. Because $\Sigma$ is strongly tensorial the first $n$ outputs of $\sigma$ already lie in $M$, we have only to check its rightmost output.  The quantum group case is in \cite{BrzMa:gau}. \endproof
 
As a corollary, if $\Phi'$ is a second trivialisation then the associated local section is the local gauge transformation $\gamma=\Phi'*\Phi^{-1}$ such that $\Phi'=\Phi^\gamma$, proving transitivity of the action of local gauge transformations. Next, if $\omega$ is a connection with associated projection $\Pi$, we extend $\id-\Pi$ as a left $P$-module morphism to $\Omega^nP$ in the canonical way (projecting each copy of $\Omega^1P$). We then define the covariant derivative on pseudotensorial forms as $D\Sigma =(\id-\Pi)\circ d\Sigma$. It is clear from equivariance of $\Pi$ and $d$ that $D\Sigma$ is again equivariant. It remains to see, however, when it descends to strongly tensorial forms:

\begin{propos} let $P,B$ be a braided principal bundle. Then the covariant derivative $D$ associated to connection $\omega$ sends strongly tensorial $n$-forms  to strongly tensorial $n+1$-forms {\em iff} $\omega$ is strong.  If the bundle is trivial and the connection is strong then 
\[ D(\sigma*\Phi)=(\nabla\sigma)*\Phi;\quad \nabla\sigma=d\sigma+(-1)^{n+1}\sigma * A,\]
where $A$ is the corresponding gauge field. We call $\nabla$ the {\em covariant derivative} on local sections $\sigma$.
\end{propos}
\begin{figure}
\[ \epsfbox{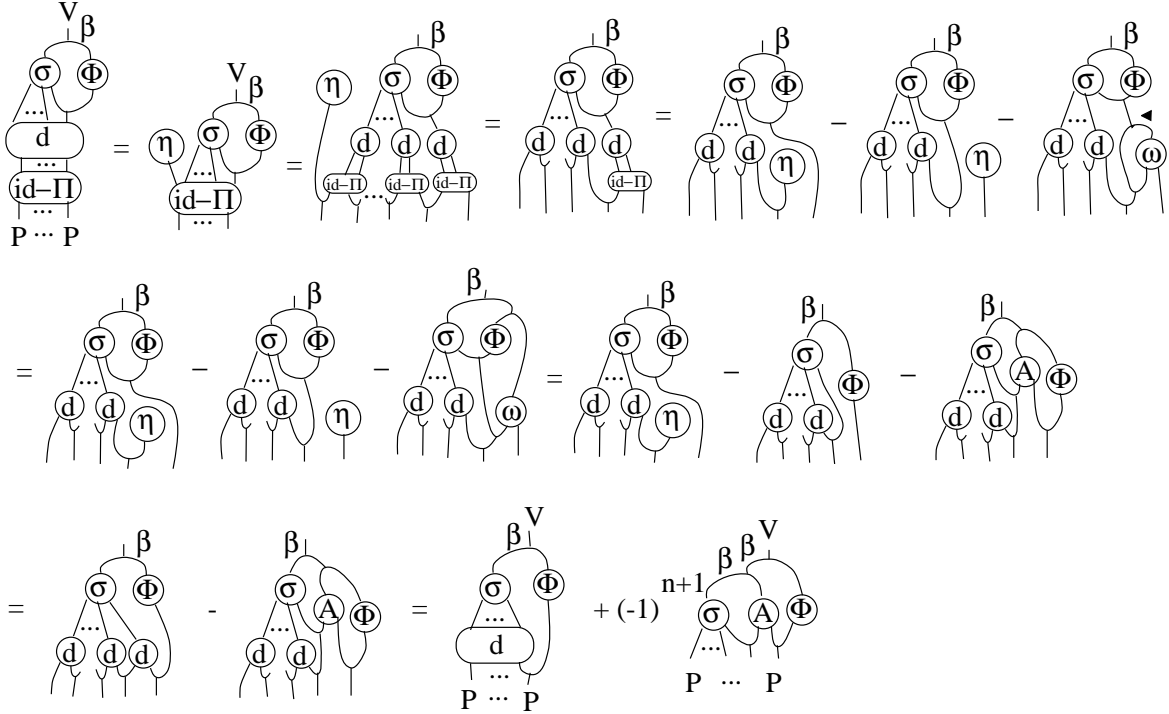}\]
\caption{Proof that $D$ for a strong connection descends to $\nabla$ on local sections.}
\end{figure}
\proof This is shown in Figure~13. We consider strongly tensorial forms of the form $\sigma*\Phi$ in the case of a trivial bundle. We use $d$ in the form on $P^{\tens n+1}$ in Figure~1. There is a signed sum with $\eta$ in all positions, but only this first terms survives the next step: we apply $\id-\Pi$ on $P^{\tens n+2}$ by mapping this to $P(\Omega^1P)^n$ via $d$, applying $\id-\Pi$ to each $\Omega^1P$ and multiplying up to return to $P^{\tens n+2}$. Moreover, $(\id-\Pi)$ is the identity on $\Omega^1M$. We then insert the form of $d$ and $(\id-\Pi)$ on the remaining output of $\sigma$ and $\Phi$, and use covariance of $\Phi$. This takes us the second line. We then insert the form of $\omega$ in terms of a gauge field $A$, in the case of a strong connection and combine the resulting first two terms as a final $d$. We finally identify the resulting first term as $(d\sigma)*\Phi$. For the second term we write $d=\eta\tens\id-\id\tens\eta$ for each $d$, but only $-\id\tens\eta$ contributes each time because $\sigma$ has its output in $\Omega^nM$. Hence we obtain the result for $D(\sigma*\Phi)$ as stated. Moreover,  
for a general bundle replace $\sigma*\Phi$ by a strongly tensorial form $\Sigma$ in the first line in Figure~13. We see that $D\Sigma$ is something in a tensor power of $M$ multiplying its rightmost factor with $(\id-\Pi)d$ acting on the rightmost output of $\Sigma$. So the result is strongly tensorial when $\omega$ is strong, just by the definition (C3). For the converse direction, take $V=P$ and consider the strongly tensorial 0-form $\id:P\to P$. We require $D(\id)=(\id-\Pi)\circ d$ to be strongly tensorial as well, i.e. to factor through $(\Omega^1M)P$. This is just the definition that $\omega$ is a strong connection. The quantum group case is basically in \cite{BrzMa:gau} with the explicit characterisation of strongness in this case in \cite{Haj:str}. \endproof

Moreover, if $\gamma$ is a gauge transformation then $\sigma^\gamma=\sigma*\gamma$ corresponds to the same $\Sigma$ when computed with respect to the trivialisation $\Phi$ as $\sigma$ provides when computed with the gauge transformed $\Phi^\gamma$.
One may confirm that $\nabla$ is then indeed covariant when $\sigma,A$ transform as here and in Section~3.  It is also easy to see that 
\[ D^2(\sigma*\Phi)=(\nabla^2\sigma)*\Phi=-(\sigma*F)*\Phi;\quad F=dA+A*A,\]
where the {\em curvature form} $F:V\to \Omega^2M$ obeys the Bianchi identities $dF+A*F-F*A$ and $F(A^\gamma)=\gamma^{-1}*F(A)*\gamma$. These local computations  follow just the same form as given in detail in \cite{BrzMa:gau}. 

To complete our picture we would like to think of the morphisms  $\sigma$ 
as the local form of cross sections of an associated fiber bundle.  Following the quantum group case in \cite{BrzMa:gau}, we define the associated bundle as a `fixed subspace' $E$ of $P\tens V$. In the quantum group case \cite{BrzMa:gau} this was done in a way that avoids any kind of commutativity of the quantum group function algebra. Some version of this should also be natural in the braided setting. 
On the other hand,  for braided groups, there is in fact a natural `commutativity condition' which can be imposed\cite{Ma:bg} and which holds for many examples. We develop this version now, as a complement to a setting more in the line of \cite{BrzMa:gau}. Recall from \cite{Ma:bg} that $B$ is {\em braided commutative} with respect to a right $B$-comodule $V$ (say) if
\[ \epsfbox{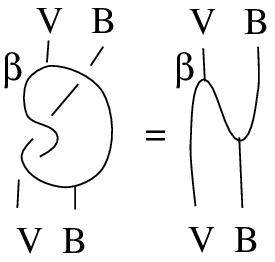}\]
holds. Here $\beta$ is the coaction of $B$ on $V$. The same diagram applies up-side-down for braided cocommutativity with respect to a module\cite{Ma:bra}.

Following essentially the proof in \cite{Ma:tra} (turned up-side-down) that the category of comodules with respect to which $B$ is braided-commutative is closed under tensor product, one finds:

\begin{propos} Let $P,V$ be braided $B$-comodule algebras and $B$ braided-commutative with respect to $V$. Then the braided tensor product coaction makes the braided tensor product algebra $P\und\tens V$ into a braided $B$-comodule algebra. 
\end{propos}
\begin{figure}
\[\epsfbox{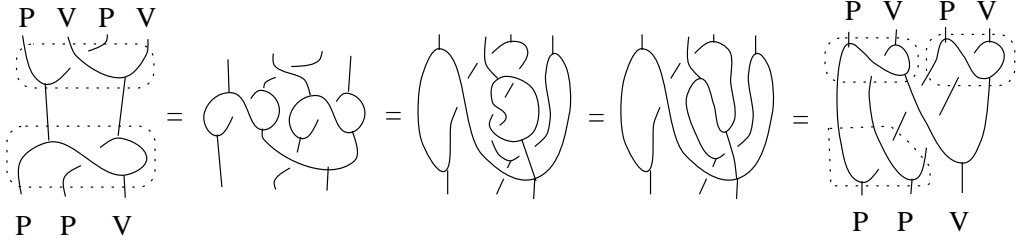}\]
\caption{Proof that braided tensor product $P\und\tens V$ becomes a braided comodule algebra when $B$ is braided-commutative with respect to $V$}
\end{figure}
\proof This is given in Figure~14, with all $\epsfbox{deltafrag.eps}$ nodes denoting the coaction.
We use the comodule homomorphism property for each of $P,V$, deform the diagram to the required form and  apply the above cocommutativity condition. We obtain the comodule homomorphism property for $P\und\tens V$ with the braided tensor product algebra structure (upper left box) and the braided tensor product coaction (lower left box). \endproof

Motivated by this, we define $E=(P\und\tens V)^B$ to be the `fixed point' equaliser object {\em  whether or not}  $B$ is braided-commutative with respect to $V$ or  $V$ a braided comodule algebra. We require only that $V$ comodule equipped with a morphism $\eta_V:\id\to V$ fixed under the coaction. Then we still have an object $E$ which remains a left $M$-module, though no geometrical picture of it as the `coordinate ring' of the total space of the bundle. It is clear that our suppressed morphism $M\to P$ induces a morphism $M\to E$, which we also suppress. We call $E$ the braided fiber bundle associated to braided principal bundle $P,B$ and the comodule $V$. 
We also define a {\em cross section} of $E$ to be a left $M$-module morphism $\vecs:E\to M$ such $\vecs\circ\eta_E=\eta_M$. 

\begin{propos} Let $P,B$ be a braided principal bundle and $E$ associated to it with fiber $V$. Suppose that the braided antipode $S$ of $B$ is invertible. Then pseudotensorial 0-forms $\Sigma$ such that $\Sigma\circ\eta_V=\eta_P$ are in 1-1 correspondence with cross sections $\vecs$, via
\[ \epsfbox{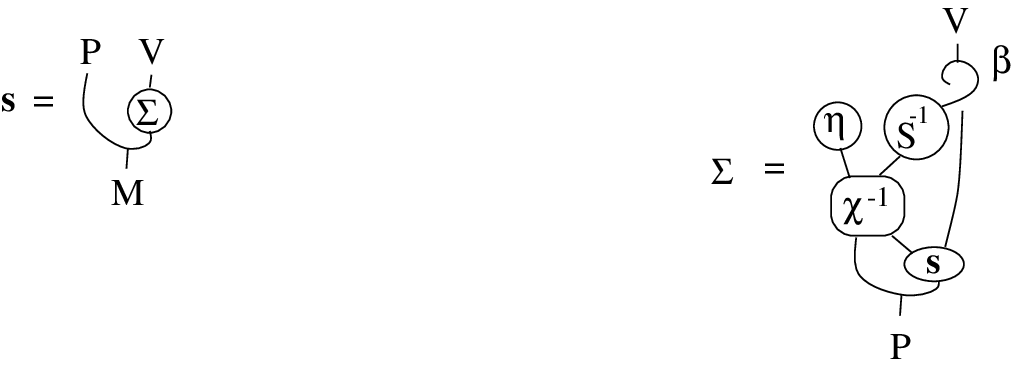}\]
\end{propos}
\begin{figure}
\[ \epsfbox{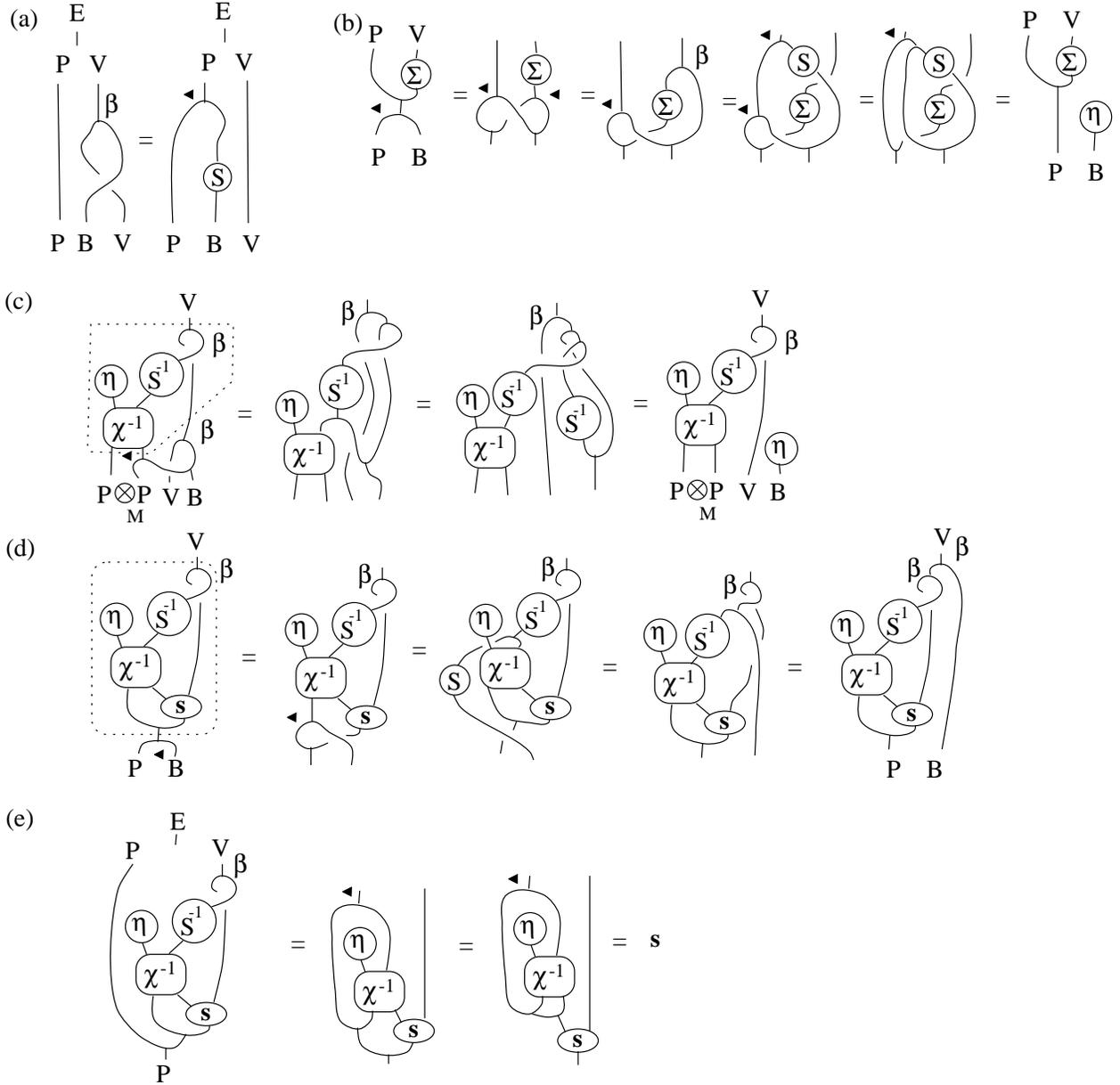}\]
\caption{(a) Identity in $E$ used in (b) construction of cross section from pseudotensorial 0-form
and (c)-(e) converse construction}
\end{figure}
\proof This is shown in Figure~15. Part (a) begins with a useful identity for coaction on $E$. It is a rearrangement of the condition that $E$ is fixed under the braided tensor product coaction on $P\tens V$. Part (b) shows that $\vecs$ defined by $\Sigma$ factors through $M$. We use the comodule homomorphism property, that $\Sigma$ is pseudotensorial and part (a). We then use the comodule property and cancel the resulting antipode loop. The other properties of $\vecs$ are immediate. Part (c) checks that the combination involving $\chi^{-1}$ used in the construction of $\Sigma$ from $\vecs$ has its output in $E$, so that the stated formula for $\Sigma$ makes sense. The formula also depends on $\vecs$ being an $M$-module morphism to make sense as stated, since the output of $\chi^{-1}$ is in $P\tens_MP$. We apply the tensor product coaction to $P\tens V$, use the comodule property and covariance of $\chi^{-1}$ from Figure~3(b). We then use that $S^{-1}$ is a braided anticoalgebra morphism\cite{Ma:introp} and cancel the resulting $S^{-1}$ twisted loop. Part (d) verifies that $\Sigma$ defined from $\vecs$ is indeed pseudotensorial. We use the covariance of $\chi^{-1}$ from Figure~3(b), the braided anticoalgebra morphism property of $S^{-1}$ and finally the comodule property. That these constructions are inverse is immediate from one side. From the other side, starting with $\vecs$ we define a pseudotensorial form and then a cross section from it in part (e). We use part (a) and to recognise the combination in Figure~3(c) to recover $\vecs$. The correspondence for the quantum group case is in \cite{BrzMa:gau} in one direction and the converse explicitly in \cite{Brz:tra}. \endproof

\begin{corol} Let $P,B$ be a braided principal bundle and set $V=B_R$, the coregular representation. Then the associated fiber bundle $E=(P\und\tens B_R)^B\isom P$ as braided algebras. Moreover,  $P$ is a trivial braided principal bundle {\em iff} $P$ as an associated fiber bundle admits a cross section $P\to M$ such that the corresponding $B_R\to P$ is convolution invertible.
\end{corol}
\begin{figure}
\[ \epsfbox{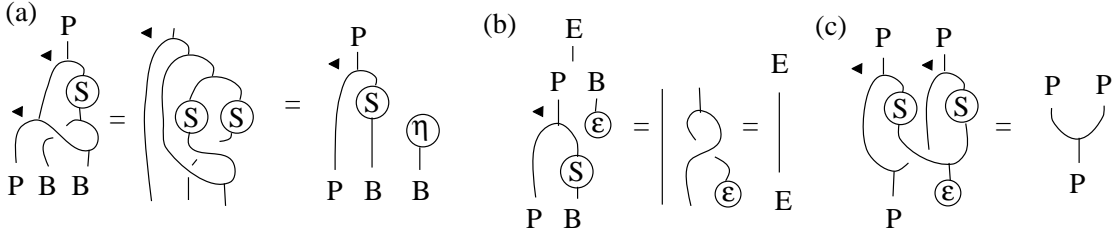}\]
\caption{Proof that the associated fiber bundle $(P\und\tens B_R)^B$ to the right coregular representation is isomorphic to $P$}
\end{figure}
\proof This is shown in Figure~16. We consider $(\id\tens S)\circ\cora:P\to P\tens B$ and show in part (a) that this indeed factors through $P\to E$. We just use the braided antimultiplicativity of $S$ and cancel the resulting antipode loop. The inverse map $E\to P$ is just $\id\tens\eps$. That this is the inverse on one side is immediate. The inverse on the other side is in part (b). We use Figure~15(a). Finally, in part (c) we apply the isomorphism, the product in $P\und\tens B$ and the inverse isomorphism and recover the product of $P$.  We then use the preceding proposition. \endproof

Moreover, we would expect that fiber bundles associated to trivial principal bundles should be trivial as well: 

\begin{propos} If $P$ is trivial with trivialisation $\Phi$ then $E\isom M\tens V$ as objects in the category, via the morphisms
\[ \epsfbox{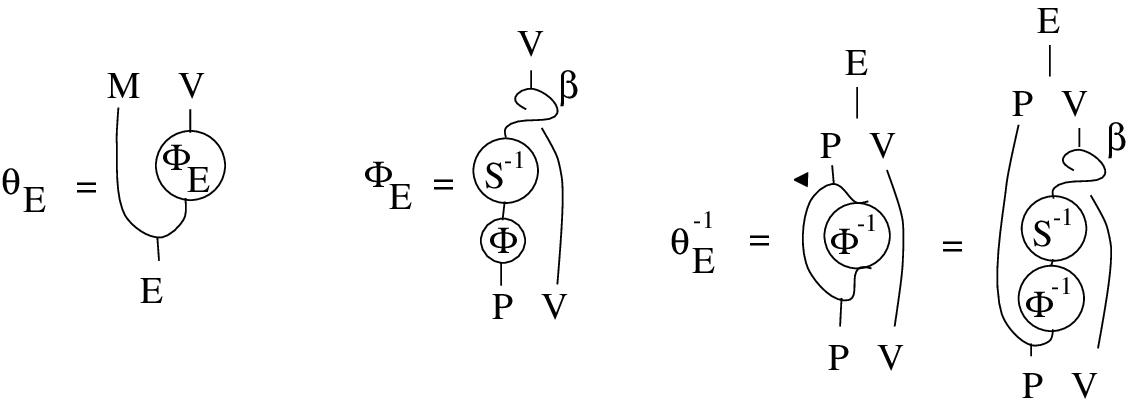}\]
We say that $E$ is a trivial   associated braided fiber bundle with trivialisation $\Phi_E:V\to E$. 
\end{propos}
\begin{figure}
\[ \epsfbox{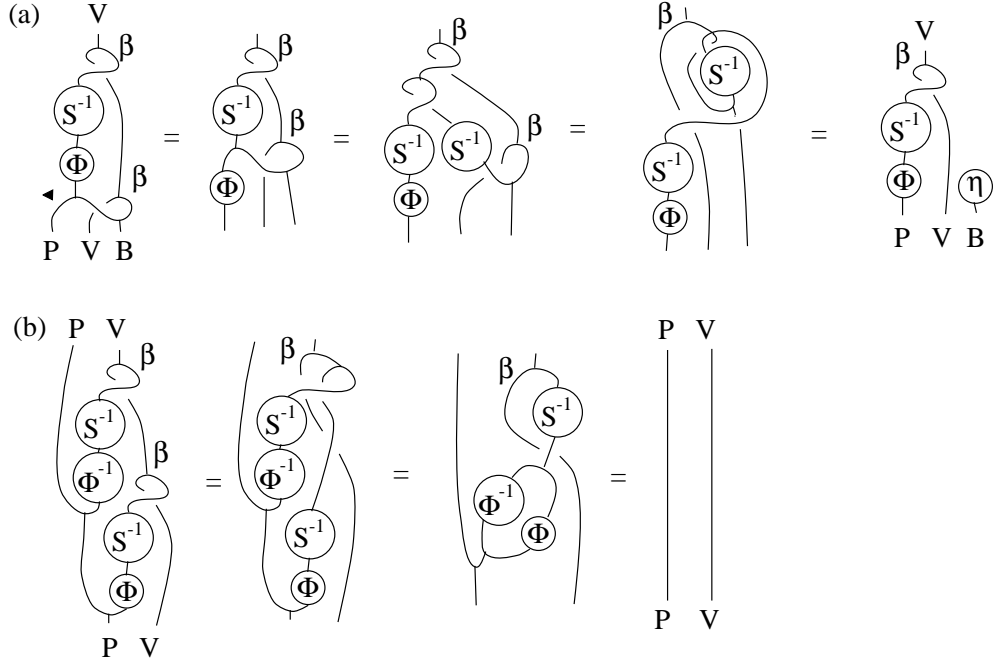}\]
\caption{Construction of trivialisation of fiber bundle associated to a trivial bundle, (a) proof that $\Phi_E:V\to E$ (b) that $\theta_E,\theta_E^{-1}$ are inverse}
\end{figure}
\proof This is given in Figure~17. Part (a) shows that $\Phi_E$ factors as claimed through $E$. We apply the braided tensor product coaction, use the covariance of $\Phi$, the braided anticomultiplicativity of $S^{-1}$ from \cite{Ma:introp} and the comodule property. We then cancel the resulting $S^{-1}$ twisted loop. Part (b) checks that $\theta_E$ and $\theta_E^{-1}$ are inverse. We use the comodule property and braided anticomultiplicativity, allowing us to cancel the resulting $\Phi^{-1},\Phi$ loop. The proof from the other side is the same with $\Phi,\Phi^{-1}$ interchanged. The quantum group case is in \cite{BrzMa:gau}. \endproof

We are then in a position to understand the above morphisms $\sigma$ as sections of the corresponding associated vector bundle. 

\begin{corol} Let $E$ be a trivial associated braided fiber bundle as above. Cross sections $\vecs:E\to M$ are in 1-1 correspondence with local sections $\sigma:V\to M$ such that $\sigma\circ\eta=\eta_M$, via
\[ \epsfbox{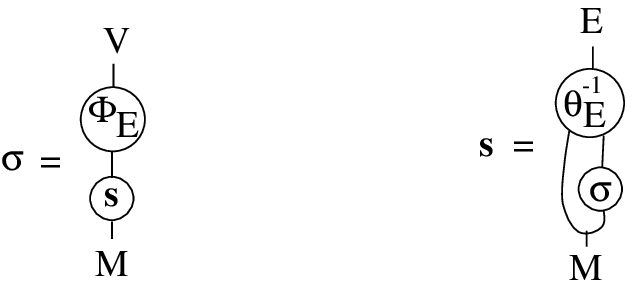}\]
\end{corol}
\proof This follows in principle from Propositions~4.1 and~4.4. For a direct proof, let $\vecs$ be defined from $\sigma$ as stated. It is clearly a left $M$-module morphism due to associativity of the product in $M$, and its restriction to $M$ is the identity since $\sigma\circ\eta=\eta$. Moreover, we recover $\sigma$ from $\vecs$ as $\vecs\circ \Phi_E =\cdot\circ(\id\tens\sigma)\circ\Phi_E=\sigma$ since $\theta_E\circ(\eta\tens\id)=\Phi_E$. Conversely, let $\sigma$ be defined from $\vecs$ as stated. Then $\cdot\circ(\id\tens\vecs)\circ(\id\tens\Phi_E)\circ\theta_E^{-1}=\vecs$ because $\vecs$ is a left $M$-module morphism, so we can move the product to $\Phi_E$ where it gives $\theta_E$. Note that we are not limited to scaler cross sections here or in Proposition~4.4; one can consider also $\vecs:E\to \Omega^nM$ corresponding in a similar way with suitable local sections $\sigma:V\to\Omega^nM$. The quantum group case is in \cite{BrzMa:gau}. 
\endproof

When $P$ is a braided tensor product principal bundle as in  Example~3.2, 
an associated braided fiber bundle also has a tensor product form:

\begin{example} Let $P=M\und\tens B$ as in Example~3.2 and $V$ a braided $B$-comodule algebra. Then $E\isom M\und\tens V$  the braided tensor product algebra.
\end{example}
\begin{figure}
\[ \epsfbox{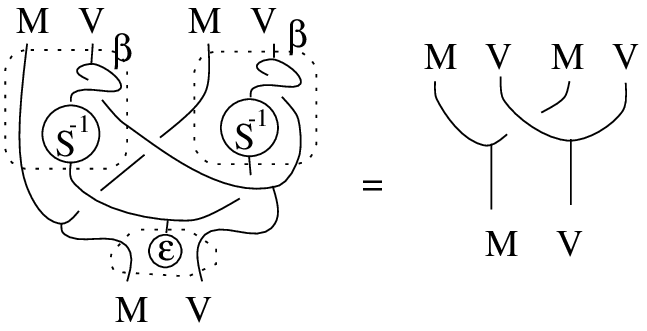}\]
\caption{Fiber bundle associated to braided tensor product bundle $M\und\tens B$ is $M\und\tens V$}
\end{figure}
\proof The proof is shown in Figure~18. We compute the product on $M\tens V$ induced by the isomorphism $\theta_E$ in Proposition~4.4. Here $\Phi=\id$. The upper boxes are $\theta_E$. We then make the product in $P\und\tens V=(M\und\tens B)\und\tens V=M\und\tens (B\und\tens V)$ and apply $\theta_E^{-1}$ in the first form in Proposition~4.4. Here $\Phi^{-1}=S$ and $\cora=
\id\tens\Delta$ so that $\theta_E^{-1}$ collapses to $\eps$. In fact, $M$ does not really enter here; one has also $(B\und\tens V)^B\isom V$ via such an isomorphism. \endproof

Finally, we return to our geometrical picture made possible by the braided theory (in contrast to the quantum group gauge theory in \cite{BrzMa:gau}. Indeed, braided-comodule algebras $V$ with respect to which $B$ is braided-commutative are not uncommon. For example, for any $B$ obtained\cite{Ma:bg} by transmutation $B=B(A,A)$ from a dual quasitriangular Hopf algebra $A$, one knows that any $A$-comodule algebra $V$ becomes via the same linear map (here $B=A$ as coalgebras) 
a braided $B$-comodule algebra. This is part of the categorical definition of transmutation as inducing a monoidal functor $A$-comodules to braided $B$-comodules\cite{Ma:bg}. Thus, although the quantum adjoint coaction is {\em not} in general a comodule algebra structure (for a noncommutative Hopf algebra) it becomes  after transmutation the braided adjoint coaction which, due to braided commutativity, is\cite{Ma:lie} a braided comodule algebra structure. (The proof of this in \cite{Ma:lie} is for the adjoint action and should be turned up-side-down to read for the braided adjoint coaction). Hence we have a fully geometrical picture of the adjoint bundle $E=(P\und\tens B_{\Ad})^B$ as a braided fixed point algebra.

\section{Example: Anyonic Gauge Theory}

Here we study the local theory for what is probably the simplest truly
braided case, namely in the braided category of $\Z_3$-graded or 3-anyonic
vector spaces\cite{Ma:any}. Objects are $\Z_3$-graded vector spaces and the braiding
is
\[ \Psi(v\tens w)=q^{|v||w|}w\tens v\]
where $q^3=1$ and $v,w$ are homogeneous of degree $|\ |$. We will study in detail the simplest case where $M$ and $B$ are 1-dimensional, i.e. something like an anyonic line bundle over an anyonic line. Of course, many other examples of the theory are equally possible, including $q$-deformations of the usual geometrical setting. Such examples, and non-trivial global bundles involving them, will be presented elsewhere. 

We work over a ground field $k$ of characteristic 0 and suppose that we have $q\in k^*$ such that $q^3=1$ and $q,q^2\ne 1$. Then let $M=k[\theta]/\theta^3$, which is the anyonic line from \cite{Ma:any}. The
degree of $\theta$ is 1. We also take $B=k[\xi]/\xi^3$ as another copy of the anyonic line, as a braided group with $\Delta \xi=\xi\tens 1+1\tens\xi$ and $\eps \xi=0$, $S\xi=-\xi$. The braided coproduct homomorphism property means
\[ \Delta(\xi^2)=(\xi\tens 1+1\tens\xi)(\xi\tens 1+1\tens\xi)=\xi^2\tens 1+(1+q)\xi\tens\xi+1\tens\xi^2\]
where we multiply in the braided tensor product algebra. We let $P=M\und\tens B$ the braided tensor product bundle. This is the algebra generated by $\xi,\theta$ with the relations
\[ \theta^3=0,\quad \xi^3=0,\quad  \xi\theta=\cdot\Psi(\xi\tens\theta)=q\theta\xi \]
and forms a $B$-comodule algebra with $\theta\mapsto \theta\tens 1$ and $\xi\mapsto \xi\tens 1+ 1\tens \xi$.

We compute first the universal differential forms on the base. Thus, 
\align{\Omega^1M&&\equad =\span\{\theta\tens\theta^2, \theta^2\tens\theta| 1\tens\theta-\theta\tens 1, \theta^2\tens\theta^2| 1\tens\theta^2-\theta^2\tens 1, \theta\tens\theta-\theta^2\tens 1\}\\
&&=\span\{ \theta d\theta^2, \theta^2d\theta| d\theta, \theta^2d\theta^2| d\theta^2, \theta d\theta\}}
in degrees $0|1|2$. This is a 6-dimensional vector space.  We give a basis for the elements of $M\tens M$ which are in the kernel of the product map, or equivalently (the second list) as elements of $\Omega^1M$ spanned by $MdM$. 

Then an anyonic gauge field means a degree-preserving map (i.e. morphism) $A:B\to\Omega^1M$ such that $A(1)=0$. We see that the  space $\CA_s$ of local gauge fields (or strong connections) is 4-dimensional, being governed by
4 parameters $a_i,b_i\in k$:
\[ A(1)=0,\quad A(\xi)=a_1d\theta+ a_2\theta^2d\theta^2,\quad A(\xi^2)=b_1d\theta^2+ b_2 \theta d\theta.\]
The curvature 2-form $F=dA+A*A$ comes out as 
\align{&&\quad F(1)=0,\quad F(\xi)=d A(\xi)=a_2 d\theta^2 d\theta^2 \\
F(\xi^2)\equad &&=dA(\xi^2)+(1+q)A(\xi)A(\xi)=(b_2+(1+q)a_1^2)d\theta d\theta+(1+q)a_1a_2(\theta^2d\theta^2 d\theta+(d\theta)\theta^2d\theta^2)\\
&&=(b_2+(1+q)a_1^2)d\theta d\theta + (1+q)a_1a_2\theta^2(d\theta^2 d\theta-d\theta d\theta^2)}
One may verify the Bianchi identity $dF+A*F-F*A$ as a useful check of the calculation. Moreover, the space of flat connections is parametrized by 
\[ a_1,b_1; \quad a_2=0,\quad b_2=-(1+q)a_1^2\]
i.e. 2-dimensional. 

A local gauge transformation means a convolution-invertible morphism $\gamma:B\to M$ such that $\gamma\circ\eta=\eta$. In fact, convolution-invertibility is automatic in the present case. So the group $\CG$ is 2-dimensional, being governed by 2 parameters $c_i\in k$:
\[ \gamma(1),\quad \gamma(\xi)=c_1\theta,\quad \gamma(\xi^2)=c_2\theta^2.\]
The group composition law is determined from
  \align{&&\equad  \gamma*\gamma'(1)=1,\quad \gamma*\gamma'(\xi)=(c_1+c_1')\theta,\quad \\
 &&\gamma*\gamma'(\xi^2)=\gamma(\xi^2)+\gamma'(\xi^2)+(1+q)\gamma(\xi)\gamma'(\xi)= (c_2+c_2'+(1+q)c_1c_1')\theta^2,}
i.e. 
\[ (c_1,c_2)+(c_1',c_2')=(c_1+c_1',c_2+c_2'+(1+q)c_1c_1')\]
which is isomorphic to the Abelian group $k^2$ under a change of coordinates to $c_2^{\rm new}=c_2-{(1+q)\over 2}c_1^2$ (our original coordinates, however, give simpler formulae below). We see that $\CG\isom k^2$ as an Abelian group.

The action of a gauge transform on a gauge field is   $A^\gamma=\gamma^{-1}*A*\gamma+\gamma^{-1}*d\gamma$, which comes out as
\align{&&\gamma^{-1}*d\gamma(1)=0,\quad \gamma^{-1}*d\gamma(\xi)=c_1 d\theta,\quad \gamma^{-1}*d\gamma(\xi^2)=c_2d\theta^2 -(1+q)c_1^2\theta d\theta\\
&&\gamma^{-1}*A*\gamma(1)=0,\quad \gamma^{-1}*A*\gamma(\xi)=A(\xi)\\
&&\gamma^{-1}*A*\gamma(\xi^2)=A(\xi^2)+(1+q)\gamma^{-1}(\xi)A(\xi)+(1+q)A(\xi)\gamma(\xi)\\
&&\qqquad\qquad =A(\xi^2)-(1+q)c_1a_1\theta d\theta+(1+q)a_1c_1(d\theta)\theta}
where we used $d1=0$ and 
\align{(\Delta\tens\id) \Delta\xi\equad &&=\xi\tens 1\tens 1+1\tens\xi\tens 1+1\tens 1\tens\xi\\
(\Delta\tens\id)\Delta\xi^2\equad&&=\xi^2\tens 1\tens 1+1\tens\xi^2\tens 1+1\tens 1\tens\xi^2+(1+q)(1\tens\xi\tens\xi+\xi\tens 1\tens\xi+1\tens\xi\tens\xi}
for the double convolution. Thus
\align{&&A^\gamma(1)=0,\quad A^\gamma(\xi)=(a_1+c_1)d\theta +a_2 \theta^2d\theta^2\\
&&A^\gamma(\xi^2)=(b_1+c_2+(1+q)a_1c_1)d\theta^2+(b_2-(1+q)(c_1^2+2c_1a_1))\theta d\theta}
i.e.
\[ \pmatrix{a_1\cr a_2\cr b_1\cr b_2}\mapsto \pmatrix{a_1+c_1\cr a_2\cr b_1+c_2+(1+q)a_1c_1\cr b_2-(1+q)(c_1^2+2c_1a_1)}\]
We see that by a gauge transformation we can set $a_1=0$ using $c_1$ and $b_1=0$ using $c_2$, leaving $a_2,b_2$ as parameters. Hence the moduli space  $\CA_s/\CG$ is 2-dimensional. The moduli space of flat connections up to gauge transformations is zero dimensional since $a_2=0$ and $b_2$ is determined from $a_1$. Indeed, we see that every flat connection is gauge equivalent to the trivial connection with zero gauge field.

Next we consider associated vector bundles and their local sections. A braided $B$-comodule means in our case a $\Z_3$-graded vector space $V$ with a degree-preserving map of the form 
\[ v\mapsto v\tens 1+ \beta(v)\tens\xi+\beta'(v)\tens\xi^2\]
The first term is dictated by the counity axiom for comodules. Moreover, the requirement of a coaction  requires us to equate the repeated coaction with the action of $\Delta$, so
\align{&&\equad v\tens 1\tens 1+\beta(v)\tens\xi\tens 1+\beta'(v)\tens\xi^2\tens 1+\beta(v)\tens 1\tens \xi+\beta^2(v)\tens\xi\tens\xi+\beta'\beta(v)\tens\xi^2\tens\xi\\
&&+\beta'(v)\tens1\tens\xi^2+\beta\beta'(v)\tens\xi\tens\xi^2+\beta'^2(v)\tens\xi^2\tens\xi^2\\
&&=v\tens 1\tens 1+\beta(v)\tens(\xi\tens 1+1\tens\xi)+\beta'(v)\tens(\xi^2\tens 1+1\tens\xi^2+(1+q)\xi\tens\xi)}
This tells us that $\beta'=\beta^2/(1+q)$, $\beta'\beta=\beta\beta'=0$. So coactions of $B$ are of the form
\[ v\mapsto v\tens 1+\beta(v)\tens\xi+{\beta^2(v)\over 1+q}\tens\xi^2 \]
for $\beta$ of degree -1 such that $\beta^3=0$. Equivalently, one can say that $B$ is dually paired with another braided group of the same form, and a coaction means an action $\beta:V\to V$ of its generator.

To be concrete, we take $V=B_R=k[\xi]/\xi^3$ the right coregular representation. So $E=P=M\und\tens B$ again. The coaction corresponds to
the operator 
\[ \beta(1)=0,\quad \beta(\xi)=0,\quad \beta(\xi^2)=(1+q)\xi\]
Scalar local sections are morphisms $\sigma:V\to M$, i.e. of the form
\[ \sigma(1)=s_0,\quad \sigma(\xi)=s_1\theta,\quad \sigma(\xi^2)=s_2\theta^2\]
The space of such local sections is 3-dimensional. For a geometrical picture where $V$ is viewed as a `coordinate ring' it is natural to fix $s_0=1$, giving a 2-dimensional affine space. The covariant derivative $\nabla\sigma=d\sigma-\sigma*A$ in the presence of a gauge field is
\align{ &&\nabla \sigma(1)=0,\quad \nabla\sigma(\xi)=d\sigma(\xi)-A(\xi)=(s_1-s_0a_1)d\theta -s_0a_2\theta^2 d\theta^2\\
&&\nabla\sigma(\xi^2)=d\sigma(\xi^2)-A(\xi^2)-(1+q)\sigma(\xi)A(\xi)=(s_2-s_0b_1)d\theta^2-(s_0b_2+(1+q)s_1a_1)\theta d\theta}
It is a nice check to compute $\nabla(\nabla\sigma)=d\nabla\sigma +(\nabla\sigma)*A$ and verify that it coincides with $-\sigma*F$.

Finally, the gauge transform of a local section is $\sigma^\gamma=\sigma*\gamma$. Thus,
\[ \sigma*\gamma(1)=s_0,\quad \sigma*\gamma(\xi)=(s_1+c_1)\theta, \quad \sigma*\gamma(\xi^2)=(s_2+s_0c_2+(1+q)s_1c_1)\theta^2\]
by the same computation as for gauge transformations. So
\[ \pmatrix{s_0\cr s_1\cr s_2}\mapsto \pmatrix{s_0\cr s_1+s_0c_1\cr s_2+s_0c_2+(1+q)s_1c_1}\]
It is a nice check of the computations to verify that $\nabla^\gamma\sigma^\gamma=(\nabla\sigma)^\gamma$, where $\nabla^\gamma$ is computed with $A^\gamma$.

This completes our description of the purely anyonic model, which is probably the simplest truly braided example of braided gauge theory. There are of course many other models that one can write down. One which is not too different for the above is to take $B=k[\xi]/\xi^3$ and before and $M=N\tens k[\theta]/\theta^3$, where $N$ is an anyonic-degree 0 and (say) commutative  algebra. So $M$ is like the coordinate ring of an `anyspace' with one anyonic dimension and the remainder bosonic. We fix a complement of  $\Omega^1k[\theta]/\theta^3$ in the tensor square, namely $\span\{1\tens1|1\tens\theta|1\tens\theta^2\}$ in degrees $0|1|2$. Then we can identify
\[ \Omega^1M=(\Omega^1N)\tens \span\{1\tens1|1\tens\theta|1\tens\theta^2\}\oplus(N\tens N)\tens \Omega^1k[\theta]/\theta^3\]
In this case gauge fields $A:B\to \Omega^1M$, gauge transformations $\gamma:B\to M$ and local sections $\sigma:B\to M$ (say) take the form
\align{ &&A(1)=0,\quad A(\xi)=A_11\tens\theta+ a_1d\theta + a_2 \theta^2d\theta^2,\quad A(\xi^2)=A_21\tens\theta^2+b_1d\theta^2+ b_2\theta d\theta\\
&&\gamma(1)=1,\quad \gamma(\xi)=c_1\theta,\quad \gamma(\xi^2)=c_2\theta^2,\quad \sigma(1)=s_0,\quad \sigma(\xi)=s_1\theta,\quad \sigma(\xi^2)=s_2\theta^2}
where $A_1,A_2\in \Omega^1N$, $a_1,a_2,b_1,b_2\in N\tens N$ and $c_1,c_2,s_0,s_1,s_2\in N$.   
 
We can then make similar computations to the above. When working with $\Omega^1M$, we let $\tens$ always denote the tensor product in $M\tens M$, omitting the tensor product in $N\tens k[\theta]/\theta^3$. Then
\align{&&  \gamma^{-1}*d\gamma(\xi)=d(c_1\theta)=(dc_1)1\tens\theta+(c_1\tens 1)d\theta\\
&&\gamma^{-1}*d\gamma(\xi^2)=d(c_2\theta^2)-(1+q)(c_1\theta\tens 1)d(c_1\theta)\\
&&\qquad\qquad=(dc_2)1\tens\theta^2+(c_2\tens 1)d\theta^2-(1+q)(c_1^2\tens 1)\theta d\theta - (1+q)c_1dc_1 \theta\tens\theta\\
&& \gamma^{-1}*A*\gamma(\xi)=A(\xi)\\
&&\gamma^{-1}*A*\gamma(\xi^2)=A(\xi^2)+(1+q)A(\xi)(1\tens c_1\theta)-(1+q)(c_1\theta\tens 1)A(\xi) \\
&&\quad=A(\xi^2)+(1+q)\left(A_1(1\tens c_1)1\tens\theta^2+a_1(1\tens c_1)(d\theta)\theta-(c_1\tens 1)A_1\theta\tens\theta-(c_1\tens 1)a_1\theta d\theta\right).}
Writing $(d\theta)\theta=d\theta^2-\theta d\theta$ and $\theta\tens\theta=1\tens\theta^2+\theta d\theta-d\theta^2$, we find the 
gauge transformation law for the components of $A$ as 
\[\pmatrix{A_1\cr A_2\cr a_1\cr a_2 \cr b_1 \cr b_2}\mapsto \pmatrix{A_1+dc_1\cr A_2+dc_2-(1+q)c_1dc_1+(1+q)(A_1\cdot c_1- c_1\cdot A_1)\cr a_1+c_1\tens 1\cr a_2\cr b_1+c_2\tens 1+(1+q)c_1dc_1+(1+q)c_1\cdot A_1+(1+q)a_1\cdot c_1\cr
b_2-(1+q)(c_1^2\tens 1+c_1dc_1)-(1+q)\left(a_1\cdot c_1 +c_1\cdot (a_1+A_1)\right)}.\]
Here $\cdot$ denotes the natural $N$-bimodule structure of $N\tens N$. The composition of gauge transformations and the transformation of local sections have the same form as before.

Likewise, we can compute the curvature $F=dA+A*A$ and hence the moduli space of zero curvature gauge fields. The most explicit way is with $\Omega^2M$ as a subset of $M\tens M\tens M$, with $d$ as in Figure~1. Then
\align{F(\xi)\equad &&=dA(\xi)=1\tens A_1(1\tens\theta)+1\tens a_1(1\tens\theta-\theta\tens 1)+1\tens a_2(\theta^2\tens\theta^2)\\
&&\quad -A_1^{13}1\tens 1\tens\theta-a_1^{13}1\tens1\tens\theta+a_1^{13}\theta\tens1\tens 1-a_2^{13}\theta^2\tens 1\tens\theta^2\\
&&\quad +A_1(1\tens\theta)\tens 1+a_1(1\tens\theta)\tens 1-a_1(\theta\tens 1)\tens 1+a_2(\theta^2\tens\theta^2)\tens 1,}
where $(\ )^{13}$ denotes placement in the 1,3 position in $M^{\tens 3}$. Combining coefficients of the obvious basis of $(k[\theta]/\theta^3)^{\tens 3}$, we find 
that $F(\xi)=0$ is equivalent to $A_1=da$ and $a_1=a\tens 1$ for some $a\in A$, and $a_2=0$. Likewise
\align{F(\xi^2)\equad&&=dA(\xi^2)+(1+q)A(\xi)A(\xi)=1\tens (b_1+A_2)(1\tens\theta^2)-1\tens (b_1+b_2)(\theta^2\tens 1)+1\tens b_2(\theta\tens\theta)\\
&&\quad -(b_1 +A_2)^{13}1\tens 1\tens\theta^2+(b_1+b_2)^{13}\theta^2\tens1\tens 1-b_2^{13}\theta\tens1\tens\theta\\
&&\quad +(b_1+A_2)(1\tens\theta^2)\tens 1-(b_1+b_2)(\theta^2\tens 1)\tens 1+b_2(\theta\tens\theta)\tens 1\\
&&\quad +(1+q)[(a_1+A_1)\cdot(a_1+A_1)1\tens\theta\tens\theta-(a_1+A_1)\cdot a_11\tens\theta^2\tens 1-a_1\cdot (a_1+A_1)\theta\tens 1\tens\theta\\
&&\qquad\qqquad +a_1\cdot a_1\theta\tens\theta\tens 1-a_1\cdot a_2\theta\tens\theta^2\tens\theta^2+a_2\cdot (a_1+A_1)\theta^2\tens\theta^2\tens\theta]}
where $\cdot$ is also used to denote the product $N^{\tens 2}\tens N^{\tens 2}\to N^{\tens 3}$ of the middle two factors. Analysing $F(\xi^2)=0$ in the same way, we find that the flat connections are of the form
\[ \pmatrix{A_1\cr A_2\cr a_1\cr a_2\cr b_1\cr b_2}=\pmatrix{da\cr db+(1+q)(da^2-ada)\cr a\tens 1\cr 0\cr b\tens 1 +(1+q)a\tens a\cr -(1+q)a\tens a},\quad a,b\in N.\]
From the above gauge transformation law, we see that every flat connection is gauge equivalent to the trivial one with zero gauge field. 

Finally, the covariant derivative on local sections is $\nabla\sigma(1)=ds_0$ and
\align{&& \nabla \sigma(\xi)=d(s_1\theta)- s_0A(\xi)=(ds_1-s_0\cdot A_1)1\tens\theta+(s_1\tens 1-s_0\cdot a_1)d\theta- s_0\cdot a_2\theta^2d\theta^2\\
&&\nabla\sigma(\xi^2)=d(s_2\theta^2)-s_0A(\xi^2)-(1+q)\sigma(\xi)A(\xi)\\
&&\qquad\quad=(ds_2-s_0\cdot A_2-(1+q)s_1\cdot A_1)1\tens\theta^2+(s_2\tens 1- s_0\cdot b_1+(1+q) s_1\cdot A_1)d\theta^2\\
&&\qquad\qquad -( s_0\cdot b_2+(1+q) s_1\cdot (a_1+A_1))\theta d\theta}
using the natural left $N$-module structure of $N\tens N$. One may verify that $\nabla$ is gauge covariant and that $\nabla^2\sigma=-\sigma*F$.

This shows how anyonic gauge theory on a composite space can appear as a novel gauge theory for a multiplet of fields on the bosonic part of the space. Moreover, it is not necessary to limit ourselves to the universal differential calculus; with a more usual commutative differential calculus on $N$,  our gauge field multiplet consists of
more usual 1-forms and some auxiliary scalar functions. One may also take a more complicated structure group $B$, for example an anyonic matrix braided group\cite{MaPla:any}. Such a theory combines the features above with the more standard features of non-Abelian gauge theory. 

\baselineskip 18pt
%\bibliographystyle{unsrt}
%\bibliography{biblio}

\end{document}